\documentclass[12pt,a4paper]{article}
 \usepackage{amsmath,amstext,amssymb,amscd,color}

\newcommand{\Xcomment}[1]{}

\newtheorem{theorem}{Theorem}

\newtheorem{corollary}[theorem]{Corollary}
\newtheorem{prop}[theorem]{Proposition}

\newenvironment{proof}{\noindent{\em Proof}.\/}%
{\hfill$\Box$\medskip}

\newcommand{\refeq}[1]{(\ref{eq:#1})}

\def\Zset{{\mathbb Z}}

\def\Lscr{\mathcal{L}}

\def\Rscr{\mathcal{R}}

\def\tone{{\tilde 1}}
\def\ttwo{{\tilde 2}}
\def\bfa{{\bf a}}
\def\bfb{{\bf b}}

\def\eps{\varepsilon}

\def\sprod{\;\hbox{\unitlength=1mm\begin{picture}(4,4)%
\put(0,0){$\triangleright$}\put(1.5,0){$\triangleleft$}%
\end{picture}}\;}

\begin{document}

 \title{$B_2$-crystals: axioms, structure, models}

 \author{Vladimir I.~Danilov
\thanks{Central Institute of Economics and
Mathematics of the RAS, 47, Nakhimovskii Prospect, 117418 Moscow, Russia;
email: danilov@cemi.rssi.ru.}
 \and
Alexander V.~Karzanov
\thanks{Central Institute of Economics and Mathematics of
the RAS, 47, Nakhimovskii Prospect, 117418 Moscow, Russia; email:
akarzanov7@gmail.com}
  \and
Gleb A.~Koshevoy
\thanks{Institute for Information Transmission Problems of
the RAS, 19, Bol'shoi Karetnyi per., 127051 Moscow, Russia; email:
koshevoyga@gmail.com. }
  }

\date{}

 \maketitle

 \begin{quote}
{\bf Abstract.} We present a list of ``local'' axioms and an explicit
combinatorial construction for the regular $B_2$-crystals (crystal
graphs of irreducible highest weight integrable modules over $U_q(sp_4)$). Also a
new combinatorial model for these crystals is developed.
 \smallskip

{\em Keywords}\,: Edge-colored graph, Doubly laced algebra, Crystal
bases of representations, Path model
\smallskip

{\em AMS Subject Classification}\, 17B37, 05C75, 05E99
  \end{quote}

\section{Introduction} \label{sec:intr}

Kashiwara~\cite{Kas-90,Kas-95} introduced the fundamental notion
of a {\em crystal} in representation theory. This is an
edge-colored directed graph in which each connected monochromatic
subgraph is a finite path and there are certain interrelations on
the lengths of such paths, described in terms of a Cartan matrix
$M$; this matrix characterizes the {\em type} of a crystal. An
important class of crystals is formed by the crystals of
representations, or {\em regular} crystals; these are associated
to irreducible highest weight integrable modules (representations) over the
quantum enveloping algebra related to $M$. There are several
global models to characterize the regular crystals for a variety
of types; e.g., via generalized Young tableaux~\cite{KN-94},
Lusztig's canonical bases~\cite{Lusztig}, Littelmann's path
model~\cite{Lit-95,Littl}.

An important fact, due to Kang et al.~\cite{KKM-92}, is that a crystal
of type $M$ is regular if and only if each (maximal connected)
2-colored subgraph in it is regular (concerning the corresponding
$2\times 2$ submatrix of $M$). This stimulates a proper study of
2-colored regular crystals.

Stembridge~\cite{Stem} pointed out a list of ``local'' graph-theoretic
axioms characterizing the regular {\em simply laced} crystals. The
2-colored subcrystals of these crystals have type $A_1\times
A_1=\binom{2\;\;0}{0\;\;2}$ or type $A_2=\binom{\;\;2\;-1}{-1\;\;2}$. A
crystal of type $A_1\times A_1$ is quite simple: it is the Cartesian
product of two paths. A combinatorial analysis of regular
$A_2$-crystals is given in~\cite{A-2}; following short terminology
there, we call such crystals {\em RA2-graphs}. It is shown
in~\cite{A-2} that any RA2-graph can be obtained from an $(A_1\times
A_1)$-crystal by replacing each monochromatic path of the latter by a
so-called {\em A-sail}, which is viewed as a triangular part of a
2-dimensional square grid.

A more complicated class is formed by the regular {\em doubly
laced} crystals. In this case the 2-colored subcrystals have type
$A_1\times A_1$ or $A_2$ or $B_2=\binom{\;\;2\;-2}{-1\;\;2}$ or $C_2=\binom{\;\;2\;-1}{-2\;\;2}$.
(A regular $B_2$- or $C_2$-crystal is associated to an irreducible
highest weight integrable module over $U_q(sp_4\simeq so_5)$.) At the end
of~\cite{Stem}, Stembridge raised the problem of characterizing
the regular $B_2$-crystals in ``local'' terms and conjectured a
complete list of possible operator relations in these crystals.
This conjecture was affirmatively answered by
Sternberg~\cite{Stern}. However, no ``local'' characterizations
for regular $B_2$-crystals have been found so far.

In this paper we attempt to give an exhaustive combinatorial analysis
of regular $B_2$-crystals. There are three groups of results that we
present. First, we give an explicit combinatorial construction for a
class of 2-edge-colored graphs, which we call {\em S-graphs}. This
construction has an analogy with the above-mentioned transformation of
an ($A_1\times A_1$)-crystal into a RA2-graph. Now we use as a base just a
RA2-graph in which certain vertices are specified, called a
``decorated'' RA2-graph. Then an S-graph is obtained from a
``decorated'' RA2-graph by replacing, in a certain way, each
monochromatic path of the latter by a so-called {\em B-sail}. Such a
sail is also a part of a square grid but in general has another
shape than an A-sail.

Second, we characterize the S-graphs by ``local'' axioms. Here we mean by the ``locality'' a way of defining a graph only via requirements on the structure of small neighborhoods of vertices (in the sense that the radius and size of a neighborhood are bounded by a constant).

Third, we develop a combinatorial {\em worm model} and show that the objects
({\em worm graphs}) generated by this model are isomorphic to S-graphs.
Moreover, a nice graphical representation of these objects enables us
to prove that the finite worm graphs satisfy the conditions in
Littelmann's path model for regular $B_2$-crystals. 
As a result, we obtain that the set of finite S-graphs is just the set of regular
$B_2$-crystals, and that these crystals are characterized by the above-mentioned
``local'' axioms. (Note that, in fact, the worm model itself leads to an alternative ``local'' axiomatics for $B_2$-crystals.)

It should be noted that our construction can produce infinite graphs as well; they can be viewed as infinite analogs of regular
$B_2$-crystals (in spirit of infinite analogs of regular $A_2$-crystals
introduced in~\cite{A-2}). However, in order to characterize such graphs axiomatically, we need to add ``non-local'' axioms.

The paper is organized as follows. In Section~\ref{sec:constr} we
briefly review a structural result on RA2-graphs from~\cite{A-2},
describe the construction of S-graphs and expose some properties of
these graphs. Section~\ref{sec:axiom} gives ``local'' axioms on a
2-edge-colored graph $G$ and proves that these axioms define exactly 
the set of S-graphs. Note that a part of axioms is stated
in terms of $G$, whereas the other axioms concern the ``decorated''
RA2-graph derived from $G$. The worm model is described in
Section~\ref{sec:worm} and we prove there that the worm graphs satisfy
the above axioms and, conversely, that any S-graph can be realized as a
worm graph. An equivalence between the worm model and Littelmann's path
model for $B_2$-crystals is proved in Appendix~1. Finally, in
Appendix~2 we explain how to transform the axioms formulated for the
``decorated'' RA2-graph derived from $G$ into ``local'' axioms directly
for $G$.

This paper is self-contained, up to appealing to Littelmann's path
model and to a structural result in~\cite{A-2}.

\section{An explicit construction}
\label{sec:constr}

In this section we present an explicit combinatorial construction
producing a class of 2-edge-colored directed graphs; we call them {\em
S-graphs} (abbreviating ``sail-graphs''). We will show later that the
finite S-graphs are precisely the regular $B_2$-crystals. Each finite
S-graph is created in a certain way from a regular $A_2$-crystal (and
there is a one-to-one correspondence between these), like the latter
can be created from an $(A_1\times A_1)$-crystal, the simplest type in
the 2-colored crystals hierarchy. Throughout in the pictures we
illustrate edges of the first color by horizontal arrows directed to
the right, and edges of the second color by vertical arrows directed
up. To avoid a possible mess, the edge colors of crystals of different types
will be denoted differently.

We start with reviewing the construction of regular $A_2$-crystal
of~\cite{A-2}, describing it in a slightly different, but equivalent,
form. For $n\in \Zset_+$, let $P_n$ denote the (directed) path of
length $n$, i.e., having $n$ edges.

A regular $A_2$-crystal is determined by parameters $a,b\in\Zset_+$,
and we denote it by $C(a,b)$. To form it, take the Cartesian product
$P_a\times P_b$, or the directed 2-dimensional rectangular grid
$\Gamma=\Gamma(a,b)$ of size $a\times b$. The latter is regarded as the
2-colored digraph whose vertices correspond to the pairs $(i,j)$,
$i=0,\ldots,a$, $j=0,\ldots,b$, the edges of the first color, say,
color $\alpha$, correspond to the pairs of the form $((i,j),(i+1,j))$,
and the edges of the second color, $\beta$ say, correspond to the pairs
of the form $((i,j),(i,j+1))$. This $\Gamma$ is a regular $(A_1\times
A_1)$-crystal,
and its vertex set $V_\Gamma$ constitutes the set of {\em principal
vertices} of $C(a,b)$ (in~\cite{A-2} they are called {\em critical}
ones).

Now $C(a,b)$ is obtained by sticking to $V_\Gamma$ copies of special
2-edge-colored graphs, so-called $A$-sails. The {\em right $A$-sail} of
size $a$ is the triangular south-east part $R$ of the $a\times a$ grid
whose vertices are the integer points $(i,j)$ for $0\le j\le i\le a$,
and the {\em left $A$-sail} of size $b$ is the triangular north-west
part $L$ of the $b\times b$ grid whose vertices are the integer points
$(i,j)$ for $0\le i\le j\le b$. Both $R$ and $L$ are 2-colored digraphs
in which the edges of the first color, say, color ${\rm I}$, are all
possible pairs of the form $((i,j),(i+1,j))$, and the edges of the
second color, ${\rm II}$ say, the pairs of the form $((i,j),(i,j+1))$.
The {\em diagonal} of $R$ (of $L$), denoted by $D(R)$ (resp. $D(L)$),
consists of the points $(i,i)$, which are ordered by increasing $i$. We
take $b+1$ copies $R_0,\ldots,R_b$ of $R$ and $a+1$ copies
$L_0,\ldots,L_a$ of $L$ and:

(i) for $j=0,\ldots,b$, replace the $j$th $\alpha$-colored path $(P_a,j)$
in $\Gamma$ by $R_j$, by identifying, for $i=0,\ldots,a$, the point
$(i,j)$ in $\Gamma$ with the $i$th point $(i,i)$ in $D(R_j)$ and then deleting
the edges of this path;

(ii) for $i=0,\ldots,a$, replace the $i$th $\beta$-colored path $(i,P_b)$
in $\Gamma$ by $L_i$, by identifying each point $(i,j)$ in $\Gamma$
with the $j$th point in $D(L_i)$ and then deleting the edges of this path.

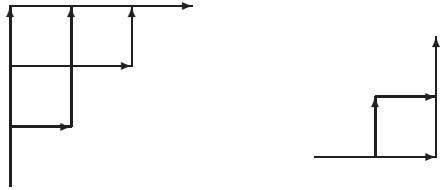
\begin{figure}[htb]
 \begin{center}
 \unitlength=.8mm \special{em:linewidth 0.4pt}
 \linethickness{0.4pt}
 \begin{picture}(100.00,35.00)(20,5)
\put(30.00,5.00){\vector(0,1){30.00}}
\put(30.00,35.00){\vector(1,0){30.00}}
\put(30.00,25.00){\vector(1,0){20.00}}
\put(30.00,15.00){\vector(1,0){10.00}}
\put(40.00,15.00){\vector(0,1){20.00}}
\put(50.00,25.00){\vector(0,1){10.00}}
\put(80.00,10.00){\vector(1,0){20.00}}
\put(100.00,10.00){\vector(0,1){20.00}}
\put(90.00,10.00){\vector(0,1){10.00}}
\put(90.00,20.00){\vector(1,0){10.00}}
 \end{picture}
 \end{center}
\caption{The left graph is $C(0,3)$ and the right graph is $C(2,0)$.}
\label{fig:fig1}
  \end{figure}

Then the resulting graph, in which the edge colors ${\rm I}$ and ${\rm
II}$ are inherited from $L$ and $R$, is just the crystal $C(a,b)$. In
\cite{A-2}, $C(a,b)$ is called the {\em diagonal-product} of $R$ and
$L$ and is denoted as $R$\sprod$L$. One can see that under this
construction $R$ and $L$ themselves are the crystals $C(a,0)$ and
$C(0,b)$, respectively (see Fig.~\ref{fig:fig1}). Also $C(a,b)$ has one {\em source}
$s$ (a zero-indegree, or minimal, vertex), which is the point $(0,0)$
in each of $\Gamma,R_0,L_0$, and one sink $t$ (a zero-outdegree, or
maximal, vertex), which is the point $(a,b)$ in $\Gamma$, $(a,a)$ in
$R_b$, and $(b,b)$ in $L_a$. The vertices of $C(a,b)$ are covered by
(inclusion-wise) maximal I-colored paths, called I-{\em strings}; these
are pairwise disjoint and each contains exactly one principal vertex.
More precisely: a vertex $(i,j)\in V_\Gamma$ belongs to the I-string
that passes the vertices $(0,j),\ldots,(j,j)$ in $L_i$ and then passes
the vertices $(i+1,i),\ldots,(a,i)$ in $R_j$ (observe that $(j,j)$ of
$L_i$ coincides with $(i,i)$ of $R_j$ and with $(i,j)$ of $\Gamma$).
Similarly, there is a natural bijection between the principal vertices
and the maximal II-colored paths, or II-{\em strings}. Note that the
parameters $a$ and $b$ are equal to the lengths of the I-string and
II-string from the source $s$, respectively (or of the II-string and
I-string to the sink $t$, respectively).

\medskip
 {\bf Finite $S$-graphs.}

 \smallskip
Next we describe the construction of the desired S-graph for parameters
$a,b\in\Zset_+$, denoted by $S(a,b)$. It is formed from $C(a,b)$ by
replacing its I-strings by right $B$-sails, and replacing
its II-strings by left $B$-sails. A $B$-sail is again a part of a
2-dimensional grid but, compared with the $A_2$ case, its structure is
somewhat more complicated and depends on the length of the I- or
II-string to be replaced by the sail as well as on the location of the
principal vertex in this string.

\begin{figure}[htb]
 \begin{center}
\unitlength=.700mm \special{em:linewidth 0.4pt}
\linethickness{0.4pt}
\begin{picture}(121.00,50)(10,10)
\put(10.00,10.00){\vector(0,1){40.00}}
\put(10.00,20.00){\vector(1,0){10.00}}
\put(20.00,20.00){\vector(0,1){30.00}}
\put(10.00,30.00){\vector(1,0){19}}
\put(30.00,30.00){\vector(0,1){20.00}}
\put(10.00,40.00){\vector(1,0){25.00}}
\put(10.00,50.00){\vector(1,0){30.00}}
\put(30.00,30.00){\circle*{2.5}} \ \ \
\bezier{16}(10.00,10.00)(20.00,20.00)(30.00,30.00)
\bezier{10}(30.00,30.00)(35.00,40.00)(40.00,50.00)
\put(55.00,35.00){\vector(1,0){10.00}}
\put(55.00,25.00){\vector(1,0){5.00}}
\put(55.00,15.00){\vector(0,1){30.00}}
\put(55.00,45.00){\vector(1,0){15.00}}
\put(55.00,15.00){\circle*{2.5}}
\put(65.00,35.00){\vector(0,1){10.00}}
\put(90.00,45.00){\vector(1,0){29}}
\put(90.00,35.00){\vector(1,0){20.00}}
\put(90.00,25.00){\vector(1,0){10.00}}
\put(90.00,15.00){\vector(0,1){30.00}}
\put(100.00,25.00){\vector(0,1){20.00}}
\put(110.00,35.00){\vector(0,1){10.00}}
\put(120.00,45.00){\circle*{2.5}}
\bezier{16}(55.00,15.00)(63.00,30.00)(70.00,45.00)
\bezier{30}(90.00,15.00)(105.00,30.00)(120.00,45.00)
\end{picture}
\end{center}
\caption{The left sails $LB(4,2)$, $LB(3,0)$, and $LB(3,3)$ (from
left to right).}
\label{fig:fig2}
  \end{figure}
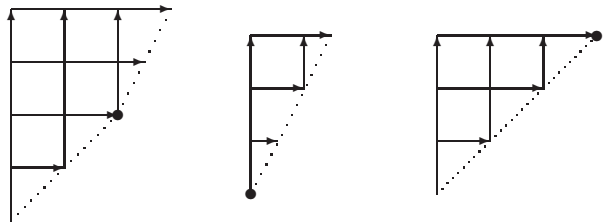

Let $r,x\in \Zset_+$ and $r\le x$. The {\em left $B$-sail} $LB(x,r)$
has the vertices identified with the integer points in the set
$\{(i,j)\colon 0\le i\le j\le x,\; j\ge 2i-r\}$ plus the half-integer
points in the set $\{(r+k+\frac 12,r+2k+1)\colon k\in \Zset_+,\; k<
\frac{x-r}{2}\}$. The point $(r,r)$ is of importance; it is called the
{\em break point} of $LB(x,r)$ (at this point the slope of the
``diagonal side'' of the sail changes from 1 to 2). The edges
correspond to all possible pairs of the form $((i,j),(i+1,j))$, to
which we assign the first color, referred to as \emph{color 1}, and all pairs of
the form $((i,j),(i,j+1))$, to which we assign \emph{color 2}. Also we
add special edges with color 1; these correspond to the pairs
$((r+k,r+2k+1),(r+k+\frac 12,r+2k+1))$ (for $k\in \Zset_+$ such that
$k< \frac{x-r}{2}$) and are called {\em left half-edges}. Instances of
left $B$-sails are illustrated in Fig.~\ref{fig:fig2} where the break points are
indicated bold.

\begin{figure}[htb]
 \begin{center}
\unitlength=.700mm \special{em:linewidth 0.4pt}
\linethickness{0.4pt}
\begin{picture}(125.00,50)(10,5)
\put(20.00,10.00){\vector(1,0){25.00}}
\put(25.00,20.00){\vector(1,0){20.00}}
\put(30.00,30.00){\vector(1,0){15.00}}
\put(35.00,40.00){\vector(1,0){10.00}}
\put(25.00,10.00){\vector(0,1){10.00}}
\put(35.00,10.00){\vector(0,1){29}}
\put(45.00,10.00){\vector(0,1){40.00}}
\put(35.00,40.00){\circle*{2.5}} \
\bezier{22}(20.00,10.00)(28.00,25.00)(35.00,40.00)
\bezier{8}(35.00,40.00)(39.00,45.00)(45.00,50.00)
\put(60.00,10.00){\vector(1,0){20.00}}
\put(80.00,10.00){\vector(0,1){39}}
\put(65.00,20.00){\vector(1,0){15.00}}
\put(70.00,30.00){\vector(1,0){10.00}}
\put(75.00,40.00){\vector(1,0){5.00}}
\put(80.00,50.00){\circle*{2.5}}
\bezier{20}(80.00,50.00)(70.00,30.00)(60.00,10.00)
\put(70.00,10.00){\vector(0,1){20.00}}
\put(95.00,15.00){\vector(1,0){30.00}}
\put(125.00,15.00){\vector(0,1){30.00}}
\put(105.00,25.00){\vector(1,0){20.00}}
\put(115.00,35.00){\vector(1,0){10.00}}
\put(105.00,15.00){\vector(0,1){10.00}}
\put(115.00,15.00){\vector(0,1){20.00}}
\put(95.00,15.00){\circle*{2.5}}
\bezier{30}(95.00,15.00)(110.00,30.00)(125.00,45.00)
\end{picture}
 \end{center}
\caption{The right sails $RB(4,1)$, $RB(4, 0)$ and $RB(3,3)$ (from
left to right).}
\label{fig:fig3}
  \end{figure}
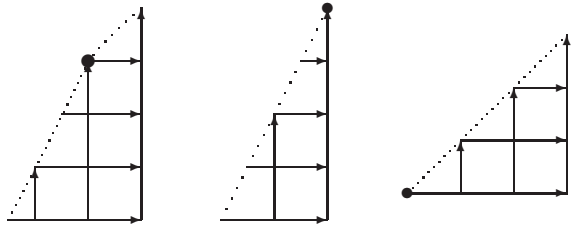

Right $B$-sails are defined symmetrically. For $q,y\in\Zset_+$ with
$q\le y$, the vertex set of the {\em right sail} $RB(y,q)$ is formed by
the integer points in $\{(i,j)\colon -y\le j\le i\le 0,\; j\le 2i+q\}$
and by the half-integer points in $\{(-q-k-\frac 12,-q-2k-1)\colon k\in
\Zset_+,\; k<\frac{y-q}{2}\}$. The break point of $RB(y,q)$ is
$-(q,q)$. The edges are the pairs of the form $((i,j),(i+1,j))$,
colored 1, the pairs of the form $((i,j),(i,j+1))$, colored 2, and the special pairs $((-q-k-\frac 12,-q-2k-1), (-q-k,-q-2k-1))$
(for $k\in \Zset_+$ such that $k<\frac{y-q}{2}$), which are colored
1 and called {\em right half-edges}. Instances of right $B$-sails
are illustrated in Fig.~\ref{fig:fig3}.

\medskip
Now we are ready to define $S(a,b)$.

For each II-string $P$ of $C(a,b)$ we proceed as follows. If $v=(i,j)$
is the principal vertex in $P$ (using the coordinates in
$\Gamma(a,b)$), then the part $P'$ of $P$ from the beginning to $v$ has
length $i$ (as $P'$ passes the vertices $(i,0),\ldots,(i,i)$ of the
right $A$-sail $R_j$), and the part $P''$ of $P$ from $v$ to the end
has length $b-j$ (as $P''$ passes the vertices $(j,j),\ldots,(j,b)$ of
the left $A$-sail $L_i$). We replace $P$ by the left $B$-sail
$LB=LB(x,r)$ with $x=b+i-j$ and $r=i$, by consecutively identifying the
vertices of $P$ (in their order) with the vertices of the ``diagonal
side'' of $LB$ (in their natural ordering starting from the minimal
vertex $(0,0)$) and then deleting the edges of $P$. So $v$ is
identified with the break point of $LB$.

For each I-string $P$ of $C(a,b)$, the procedure is similar. If
$v=(i,j)$ is the principal vertex in $P$, then the part $P'$ of $P$
from the beginning to $v$ has length $j$ (as $P'$ passes the vertices
$(0,j),\ldots,(j,j)$ of the left $A$-sail $L_i$), and the part $P''$ of
$P$ from $v$ to the end has length $a-i$ (as $P''$ passes the vertices
$(i,i),\ldots,(a,i)$ of the right $A$-sail $R_j$). We replace $P$ by
the right $B$-sail $RB=RB(y,q)$ with $y=a+j-i$ and $q=a-i$, by
consecutively identifying the vertices of $P$ with the vertices of the
``diagonal side'' of $RB$ (in their ordering starting from the minimal
vertex $(-\frac{y+q}{2},-y)$) and then deleting the edges of $P$.
Again, $v$ is identified with the break point of $RB$.

Finally, one can see that the half-integer points of the diagonals of
$B$-sails are
identified with precisely those vertices $v$ of $C(a,b)$ that belong to
left $A$-sails and lie at odd distance from their diagonals. When the
II-string passing such a $v$ is replaced by the corresponding left
$B$-sail, $v$ becomes incident to a left half-edge $(u,v)$, and when
the I-string passing $v$ is replaced by the corresponding right
$B$-sail, $v$ becomes incident to a right half-edge $(v,w)$. We merge
these half-edges into one edge $e_v=(u,w)$ (so the vertex $v$
vanishes), which inherits the color 1 from $(u,v),(v,w)$.

The resulting 2-colored digraph, with colors 1 and 2, is
just the graph $S(a,b)$. The vertices of $C(a,b)$ occurring in
$S(a,b)$, as well as the edges $e_v$ obtained by merging half-edges,
are called {\em central} elements of $S(a,b)$ (so there is a natural
bijection between the central elements and the vertices of $C(a,b)$);
the principal vertices (viz. the break points of $B$-sails) are most
important among these.

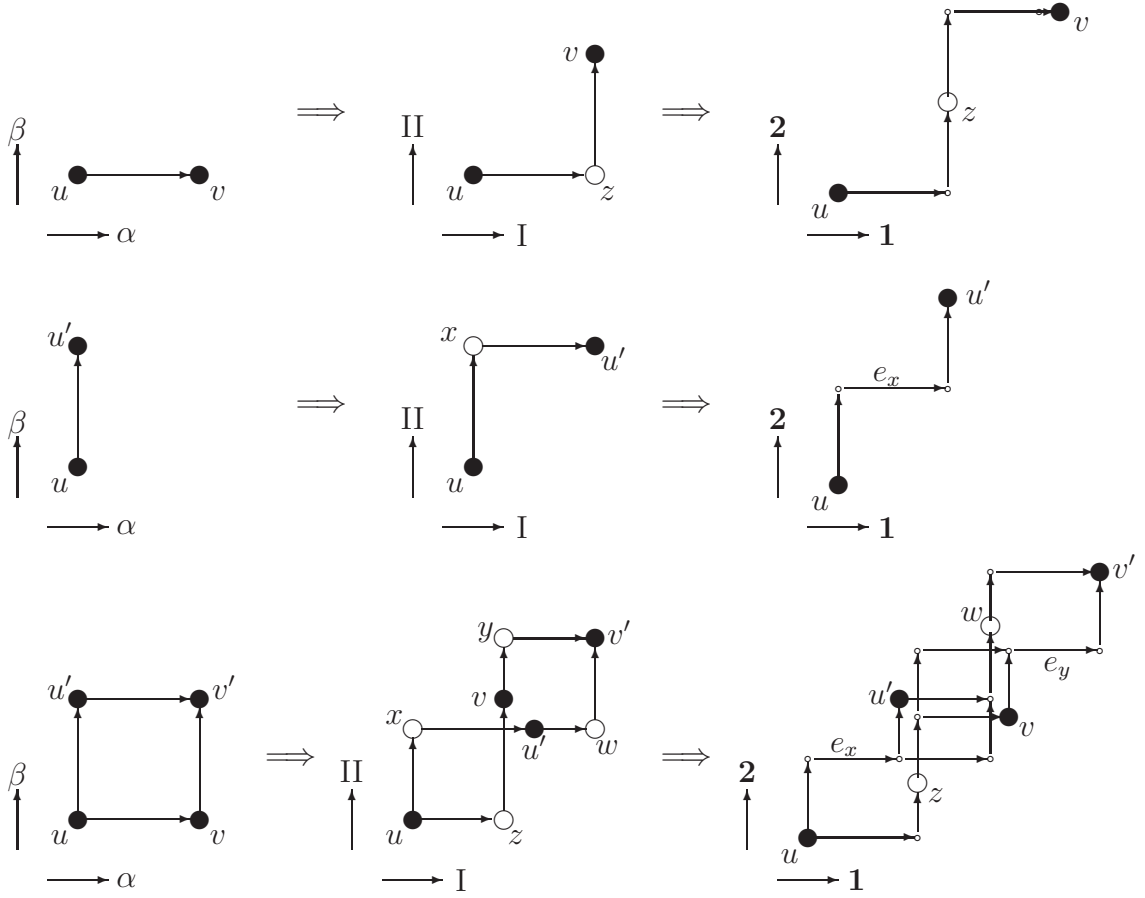
\begin{figure}[htb]
 \begin{center}
\unitlength=.800mm \special{em:linewidth 0.4pt} \linethickness{0.4pt}
\begin{picture}(180.00,40.00)(0,0)
 \put(10.00,10.00){\circle*{3.00}}
 \put(30.00,10.00){\circle*{3.00}}
 \put(10,10){\vector(1,0){18.5}}
 \put(7,7){\makebox(0,0)[cc]{$u$}}
 \put(33,7){\makebox(0,0)[cc]{$v$}}
 \put(5,0){\vector(1,0){10}}
 \put(18,0){\makebox(0,0)[cc]{$\alpha$}}
 \put(0,5){\vector(0,1){10}}
 \put(0,17){\makebox(0,0)[cc]{$\beta$}}
\put(50.00,20.00){\makebox(0,0)[cc]{$\Longrightarrow$}}
 \put(75.00,10.00){\circle*{3.00}}
 \put(72,7){\makebox(0,0)[cc]{$u$}}
 \put(95,30){\circle*{3.00}}
 \put(91,30){\makebox(0,0)[cc]{$v$}}
 \put(95,10){\circle{3.00}}
 \put(97,7){\makebox(0,0)[cc]{$z$}}
 \put(75,10){\vector(1,0){18.00}}
 \put(95,11.5){\vector(0,1){17.00}}
 \put(70,0){\vector(1,0){10}}
 \put(83,0){\makebox(0,0)[cc]{I}}
 \put(65,5){\vector(0,1){10}}
 \put(65,18){\makebox(0,0)[cc]{II}}
 \put(130,0){\vector(1,0){10}}
 \put(143,0){\makebox(0,0)[cc]{{\bf 1}}}
 \put(125,5){\vector(0,1){10}}
 \put(125,18){\makebox(0,0)[cc]{{\bf 2}}}
  \put(110.00,20.00){\makebox(0,0)[cc]{$\Longrightarrow$}}
   \begin{picture}(180.00,40.00)(-15,-5)
 \put(138.00,2.00){\circle{1.00}}
 \put(138.00,32.00){\circle{1.00}}
 \put(153.00,32.00){\circle{1.00}}
 \put(121.00,2.00){\vector(1,0){16.5}}
 \put(138.00,3.00){\vector(0,1){12.5}}
 \put(138.00,18.00){\vector(0,1){13.5}}
 \put(139.00,32.00){\vector(1,0){16.5}}
 \put(120.00,2.00){\circle*{3.00}}
 \put(117.00,-1.00){\makebox(0,0)[cc]{$u$}}
 \put(138.00,17.00){\circle{3.00}}
 \put(141.5,15.00){\makebox(0,0)[cc]{$z$}}
 \put(156.5,32.00){\circle*{3.00}}
\put(160,30.00){\makebox(0,0)[cc]{$v$}}
   \end{picture}
\end{picture}
   %

\unitlength=.800mm \special{em:linewidth 0.4pt} \linethickness{0.4pt}
\begin{picture}(180.00,48.00)(0,0)
 \put(10.00,10.00){\circle*{3.00}}
 \put(10.00,30.00){\circle*{3.00}}
 \put(10,10){\vector(0,1){18.5}}
 \put(7,7){\makebox(0,0)[cc]{$u$}}
 \put(7,32){\makebox(0,0)[cc]{$u'$}}
 \put(5,0){\vector(1,0){10}}
 \put(18,0){\makebox(0,0)[cc]{$\alpha$}}
 \put(0,5){\vector(0,1){10}}
 \put(0,17){\makebox(0,0)[cc]{$\beta$}}
\put(50.00,20.00){\makebox(0,0)[cc]{$\Longrightarrow$}}
 \put(75.00,10.00){\circle*{3.00}}
 \put(72,7){\makebox(0,0)[cc]{$u$}}
 \put(95,30){\circle*{3.00}}
 \put(98,28){\makebox(0,0)[cc]{$u'$}}
 \put(75.00,30.00){\circle{3.00}}
 \put(71,32){\makebox(0,0)[cc]{$x$}}
 \put(76.5,30){\vector(1,0){17.00}}
 \put(75,10){\vector(0,1){18.5}}
 \put(70,0){\vector(1,0){10}}
 \put(83,0){\makebox(0,0)[cc]{I}}
 \put(65,5){\vector(0,1){10}}
 \put(65,18){\makebox(0,0)[cc]{II}}
 \put(130,0){\vector(1,0){10}}
 \put(143,0){\makebox(0,0)[cc]{{\bf 1}}}
 \put(125,5){\vector(0,1){10}}
 \put(125,18){\makebox(0,0)[cc]{{\bf 2}}}
  \put(110.00,20.00){\makebox(0,0)[cc]{$\Longrightarrow$}}
   \begin{picture}(180.00,50.00)(-15,-5)
 \put(120.00,18.00){\circle{1.00}}
 \put(138.00,18.00){\circle{1.00}}
 \put(120.00,3.00){\vector(0,1){14}}
 \put(121.00,18.00){\vector(1,0){16.5}}
 \put(138.00,19.00){\vector(0,1){12.5}}
  %
 \put(120.00,2.00){\circle*{3.00}}
 \put(117.00,-1.00){\makebox(0,0)[cc]{$u$}}
 \put(138.00,33.00){\circle*{3.00}}
 \put(143.00,34.00){\makebox(0,0)[cc]{$u'$}}
 \put(128.00,20.00){\makebox(0,0)[cc]{$e_x$}}
   \end{picture}
\end{picture}
   %

\unitlength=.800mm \special{em:linewidth 0.4pt} \linethickness{0.4pt}
\begin{picture}(180.00,58.00)(0,0)
\put(10.00,10.00){\circle*{3.00}} \put(30.00,10.00){\circle*{3.00}}
\put(10.00,30.00){\circle*{3.00}} \put(30.00,30.00){\circle*{3.00}}
\put(10,10){\vector(1,0){18.5}} \put(10,30){\vector(1,0){18.5}}
\put(10,10){\vector(0,1){18.5}} \put(30,10){\vector(0,1){18.5}}
\put(7,7){\makebox(0,0)[cc]{$u$}} \put(33,7){\makebox(0,0)[cc]{$v$}}
\put(7,32){\makebox(0,0)[cc]{$u'$}}
\put(34,32){\makebox(0,0)[cc]{$v'$}}
 \put(5,0){\vector(1,0){10}}
 \put(18,0){\makebox(0,0)[cc]{$\alpha$}}
 \put(0,5){\vector(0,1){10}}
 \put(0,17){\makebox(0,0)[cc]{$\beta$}}
\put(45.00,20.00){\makebox(0,0)[cc]{$\Longrightarrow$}}
\put(65.00,10.00){\circle*{3.00}} \put(62,7){\makebox(0,0)[cc]{$u$}}
\put(85,25){\circle*{3.00}} \put(85,22){\makebox(0,0)[cc]{$u'$}}
\put(80,30){\circle*{3.00}}  \put(76,30){\makebox(0,0)[cc]{$v$}}
\put(95,40){\circle*{3.00}} \put(99,41){\makebox(0,0)[cc]{$v'$}}
\put(65.00,25.00){\circle{3.00}} \put(62,27){\makebox(0,0)[cc]{$x$}}
\put(95.00,25.00){\circle{3.00}} \put(97,22){\makebox(0,0)[cc]{$w$}}
\put(80,10){\circle{3.00}} \put(82,7){\makebox(0,0)[cc]{$z$}}
\put(80,40){\circle{3.00}} \put(77,41){\makebox(0,0)[cc]{$y$}}
\put(65,10){\vector(1,0){13.00}} \put(66.5,25){\vector(1,0){17.00}}
\put(85,25){\vector(1,0){8.5}}  \put(81.5,40){\vector(1,0){12.00}}
\put(65,10){\vector(0,1){13.5}} \put(80,11.5){\vector(0,1){17.00}}
\put(80,30){\vector(0,1){8.5}} \put(95,26.5){\vector(0,1){12}}
 \put(60,0){\vector(1,0){10}}
 \put(73,0){\makebox(0,0)[cc]{I}}
 \put(55,5){\vector(0,1){10}}
 \put(55,18){\makebox(0,0)[cc]{II}}
 \put(125,0){\vector(1,0){10}}
 \put(138,0){\makebox(0,0)[cc]{{\bf 1}}}
 \put(120,5){\vector(0,1){10}}
 \put(120,18){\makebox(0,0)[cc]{{\bf 2}}}
  \put(110.00,20.00){\makebox(0,0)[cc]{$\Longrightarrow$}}
   \begin{picture}(180.00,50.00)(-10,-5)
\put(120.00,15.00){\circle{1.00}} \put(135.00,15.00){\circle{1.00}}
\put(135.00,25.00){\circle{1.00}} \put(150.00,25.00){\circle{1.00}}
\put(150.00,46.00){\circle{1.00}} \put(138.00,2.00){\circle{1.00}}
\put(138.00,22.00){\circle{1.00}} \put(153.00,22.00){\circle{1.00}}
\put(153.00,33.00){\circle{1.00}} \put(168.00,33.00){\circle{1.00}}
\put(138.00,33.00){\circle{1.00}} \put(150.00,15.00){\circle{1.00}}
\put(120.00,3.00){\vector(0,1){11.5}}
\put(121.00,15.00){\vector(1,0){13.5}}
\put(135.00,16.00){\vector(0,1){7.5}}
\put(136.00,25.00){\vector(1,0){13.5}}
\put(150.00,26.00){\vector(0,1){10}}
\put(150.00,38.00){\vector(0,1){7.5}}
\put(151.00,46.00){\vector(1,0){16.00}}
\put(121.00,2.00){\vector(1,0){16.5}}
\put(138.00,3.00){\vector(0,1){6.5}}
\put(138.00,12.00){\vector(0,1){9.5}}
\put(139.00,22.00){\vector(1,0){13.00}}
\put(153.00,23.00){\vector(0,1){9.5}}
\put(154.00,33.00){\vector(1,0){13.5}}
\put(168.00,34.00){\vector(0,1){10.5}}
\put(138.00,23.00){\vector(0,1){9.5}}
\put(139.00,33.00){\vector(1,0){13.5}}
\put(136.00,15.00){\vector(1,0){13.5}}
\put(150.00,16.00){\vector(0,1){8.5}}
\put(120.00,2.00){\circle*{3.00}}
\put(117.00,-1.00){\makebox(0,0)[cc]{$u$}}
\put(138.00,11.00){\circle{3.00}}
\put(141.00,9.00){\makebox(0,0)[cc]{$z$}}
\put(135.00,25.00){\circle*{3.00}}
\put(132.00,26.00){\makebox(0,0)[cc]{$u'$}}
\put(168.00,46.00){\circle*{3.00}}
\put(172.00,47.00){\makebox(0,0)[cc]{$v'$}}
\put(153.00,22.00){\circle*{3.00}}
\put(156.00,20.00){\makebox(0,0)[cc]{$v$}}
\put(150.00,37.00){\circle{3.00}}
\put(147.00,39.00){\makebox(0,0)[cc]{$w$}}
\put(126.00,17.00){\makebox(0,0)[cc]{$e_x$}}
\put(161.00,30.00){\makebox(0,0)[cc]{$e_y$}}
   \end{picture}
\end{picture}
 \end{center}
\caption{In the upper line: transformation $\Gamma(1,0)
\longrightarrow C(1,0) \longrightarrow S(1,0)$. In the middle line:
transformation $\Gamma(0,1) \longrightarrow C(0,1) \longrightarrow
S(0,1)$. In the lower line: transformation $\Gamma(1,1) \longrightarrow
C(1,1) \longrightarrow S(1,1)$.}
\label{fig:fig4}
  \end{figure}

The graph $S(a,b)$ has one source $s$ and one sink $t$, which coincide
with the source and sink of $C(a,b)$, respectively. Since the I-string
and II-string beginning at $s$ are replaced by the sails $RB(a,a)$ and
$LB(b,b)$, respectively, it follows that in the graph $S(a,b)$, the
string of color 1 beginning at $s$ has length $a$ and the string of color
2 beginning at $s$ has length $b$, thus justifying the maintenance of the
parameters $a,b$. Figure~\ref{fig:fig4} illustrates three instances of the
transformation of an $(A_1 \times A_1)$-crystal into an $A_2$-crystal
and further into an S-graph.

The simplest nontrivial graphs $S(1,0)$ and $S(0,1)$ are called {\em
fundamental}. The graph $S(1,1)$ is the least S-graph where the
``big Verma relation'' is present; considering paths from the source
$s$ to the sink $t$, we observe the operator relations of degree 7
indicated in~\cite{Stem,Stern}:
 $t={\bf 1}{\bf 2}^2{\bf 1}{\bf 2}{\bf 1}{\bf 2}s=
 {\bf 1}{\bf 2}^3{\bf 1}^2{\bf 2}s=
 {\bf 2}{\bf 1}^2{\bf 2}^3{\bf 1}s=
 {\bf 2}{\bf 1}{\bf 2}{\bf 1}{\bf 2}^2{\bf 1}s$.

Figures~\ref{fig:fig4b} and~\ref{fig:fig4c} illustrate the construction of
bigger S-graphs, namely, $S(2,1)$ and $S(1,2)$, respectively. Here the
horizontal edges are directed to the right, the vertical edges are
directed up, and the central edges of the S-graphs are drawn in bold.

\begin{figure}[htb]
 \begin{center}
 \unitlength=0.8mm
\begin{picture}(180.00,160)
 \put(60,135)%
    {\begin{picture}(40,25)
   \put(10,5){\circle*{2.7}}
   \put(25,5){\circle*{2.7}}
   \put(40,5){\circle*{2.7}}
   \put(10,20){\circle*{2.7}}
   \put(25,20){\circle*{2.7}}
   \put(40,20){\circle*{2.7}}
   \put(10,5){\line(1,0){30}}
   \put(10,20){\line(1,0){30}}
   \put(10,5){\line(0,1){15}}
   \put(25,5){\line(0,1){15}}
   \put(40,5){\line(0,1){15}}
   \put(9,0.5){$a$}
   \put(24,-0.5){$b$}
   \put(39,0.5){$c$}
   \put(9,22){$A$}
   \put(24,22){$B$}
   \put(39,22){$C$}
   \put(-5,10){(a)}
     \end{picture}
    }
 \put(50,75)%
     {\begin{picture}(60,45)
   \put(10,5){\circle*{2.7}}
   \put(20,13){\circle*{2.7}}
   \put(25,20){\circle*{2.7}}
   \put(35,28){\circle*{2.7}}
   \put(40,35){\circle*{2.7}}
   \put(50,43){\circle*{2.7}}
   \put(25,5){\circle*{2}}
   \put(40,5){\circle*{2}}
   \put(10,13){\circle*{2}}
   \put(35,13){\circle*{2}}
   \put(50,13){\circle*{2}}
   \put(40,20){\circle*{2}}
   \put(25,28){\circle*{2}}
   \put(50,28){\circle*{2}}
   \put(40,43){\circle*{2}}
   \put(10,5){\line(1,0){30}}
   \put(10,13){\line(1,0){40}}
   \put(25,20){\line(1,0){15}}
   \put(25,28){\line(1,0){25}}
   \put(40,43){\line(1,0){10}}
   \put(10,5){\line(0,1){8}}
   \put(25,5){\line(0,1){23}}
   \put(35,13){\line(0,1){15}}
   \put(40,5){\line(0,1){38}}
   \put(50,13){\line(0,1){30}}
   \put(8,0.5){$a$}
   \put(21,19){$b$}
   \put(36,35){$c$}
   \put(17,15){$A$}
   \put(32,30){$B$}
   \put(52,43){$C$}
   \put(10,30){(b)}
   \put(23,0){1}
   \put(41,0){2}
   \put(7,14){3}
   \put(32,14){4}
   \put(52,10){5}
   \put(42,17){6}
   \put(22,29){7}
   \put(52,26){8}
   \put(37,44){9}
    \end{picture}
    }
 \put(20,0)%
    {\begin{picture}(100,65)
   \put(10,0){\circle*{2.7}}
   \put(18,28){\circle*{2.7}}
   \put(50,16){\circle*{2.7}}
   \put(58,44){\circle*{2.7}}
   \put(90,32){\circle*{2.7}}
   \put(98,60){\circle*{2.7}}
   \put(42,6){\circle*{2}}
   \put(74,11){\circle*{2}}
   \put(26,36){\circle*{2}}
   \put(34,49){\circle*{2}}
   \put(82,24){\circle*{2}}
   \put(66,54){\circle*{2}}
   \put(42,0){\circle*{1.2}}
   \put(74,0){\circle*{1.2}}
   \put(74,6){\circle*{1.2}}
   \put(42,16){\circle*{1.2}}
   \put(82,16){\circle*{1.2}}
   \put(10,20){\circle*{1.2}}
   \put(18,20){\circle*{1.2}}
   \put(26,20){\circle*{1.2}}
   \put(34,20){\circle*{1.2}}
   \put(74,24){\circle*{1.2}}
   \put(26,28){\circle*{1.2}}
   \put(34,28){\circle*{1.2}}
   \put(42,30){\circle*{1.2}}
   \put(50,30){\circle*{1.2}}
   \put(57,30){\circle*{1.2}}
   \put(66,30){\circle*{1.2}}
   \put(74,32){\circle*{1.2}}
   \put(82,32){\circle*{1.2}}
   \put(34,36){\circle*{1.2}}
   \put(74,40){\circle*{1.2}}
   \put(82,40){\circle*{1.2}}
   \put(90,40){\circle*{1.2}}
   \put(98,40){\circle*{1.2}}
   \put(26,44){\circle*{1.2}}
   \put(66,44){\circle*{1.2}}
   \put(34,54){\circle*{1.2}}
   \put(34,60){\circle*{1.2}}
   \put(66,60){\circle*{1.2}}
   \put(10,0){\line(1,0){64}}
   \put(42,6){\line(1,0){32}}
   \put(42,16){\line(1,0){40}}
   \put(10,20){\line(1,0){24}}
   \put(74,24){\line(1,0){8}}
   \put(18,28){\line(1,0){16}}
   \put(42,30){\line(1,0){24}}
   \put(74,32){\line(1,0){16}}
   \put(26,36){\line(1,0){8}}
   \put(74,40){\line(1,0){24}}
   \put(26,44){\line(1,0){40}}
   \put(34,54){\line(1,0){32}}
   \put(34,60){\line(1,0){64}}
   \put(10,0){\line(0,1){20}}
   \put(18,20){\line(0,1){8}}
   \put(26,20){\line(0,1){24}}
   \put(34,20){\line(0,1){40}}
   \put(42,0){\line(0,1){30}}
   \put(50,16){\line(0,1){14}}
   \put(58,30){\line(0,1){14}}
   \put(66,30){\line(0,1){30}}
   \put(74,0){\line(0,1){40}}
   \put(82,16){\line(0,1){24}}
   \put(90,32){\line(0,1){8}}
   \put(98,40){\line(0,1){20}}
 %
   \put(10,20.3){\line(1,0){8}}
   \put(10,19.7){\line(1,0){8}}
   \put(50,30.3){\line(1,0){8}}
   \put(50,29.7){\line(1,0){8}}
   \put(90,40.3){\line(1,0){8}}
   \put(90,39.7){\line(1,0){8}}
   \put(6,-3){$a$}
   \put(13,29){$A$}
   \put(51,11){$b$}
   \put(55,46){$B$}
   \put(92,30){$c$}
   \put(100,59){$C$}
   \put(38.5,4){1}
   \put(76,10){2}
   \put(22.5,34){4}
   \put(30.5,48){5}
   \put(84,23){6}
   \put(68,53){8}
   \put(13,16){$e_3$}
   \put(53,26){$e_7$}
   \put(92,42){$e_9$}
   \put(5,35){(e)}
  \end{picture}
    }
%
 \put(0,70)%
    { \begin{picture}(10,10)
   \put(0,0){\circle*{2.7}}
   \put(0,8){\line(1,0){4}}
   \put(0,0){\line(0,1){8}}
   \put(2,-2){$a$}
   \put(5,6.5){3}
  \end{picture}
    }
 \put(15,70)%
    { \begin{picture}(20,20)
   \put(8,8){\circle*{2.7}}
   \put(0,0){\circle*{2}}
   \put(0,8){\line(1,0){8}}
   \put(0,16){\line(1,0){12}}
   \put(0,0){\line(0,1){16}}
   \put(8,8){\line(0,1){8}}
   \put(2,-2){1}
   \put(10,6){$b$}
   \put(13,14){7}
  \end{picture}
    }
 \put(0,90)%
    { \begin{picture}(10,10)
   \put(8,8){\circle*{2.7}}
   \put(0,0){\circle*{2}}
   \put(0,8){\line(1,0){8}}
   \put(0,0){\line(0,1){8}}
   \put(2,-2){4}
   \put(9,6){$B$}
  \end{picture}
    }
 \put(0,110)%
    { \begin{picture}(25,25)
   \put(16,16){\circle*{2.7}}
   \put(0,0){\circle*{2}}
   \put(8,8){\circle*{2}}
   \put(0,8){\line(1,0){8}}
   \put(0,16){\line(1,0){16}}
   \put(0,24){\line(1,0){20}}
   \put(0,0){\line(0,1){24}}
   \put(8,8){\line(0,1){16}}
   \put(16,16){\line(0,1){8}}
   \put(2,-2){2}
   \put(10,7){6}
   \put(18,15){$c$}
   \put(21,22){9}
   \put(10,30){(c)}

  \end{picture}
    }
 \put(20,102)%
    { \begin{picture}(20,20)
   \put(16,16){\circle*{2.7}}
   \put(0,0){\circle*{2}}
   \put(8,8){\circle*{2}}
   \put(0,8){\line(1,0){8}}
   \put(0,16){\line(1,0){16}}
   \put(0,0){\line(0,1){16}}
   \put(8,8){\line(0,1){8}}
   \put(2,-2){5}
   \put(10,6){8}
   \put(18,14){$C$}
  \end{picture}
    }
 \put(35,95)%
    { \begin{picture}(7,7)
   \put(0,0){\circle*{2.7}}
   \put(2,0){$A$}
  \end{picture}
    }
%
%
 \put(135,100)%
    { \begin{picture}(10,10)
   \put(8,8){\circle*{2.7}}
   \put(4,0){\line(1,0){4}}
   \put(8,0){\line(0,1){8}}
   \put(3.5,8){$C$}
   \put(1,-2){9}
  \end{picture}
    }
 \put(155,60)%
    { \begin{picture}(20,20)
   \put(8,8){\circle*{2.7}}
   \put(16,16){\circle*{2}}
   \put(8,8){\line(1,0){8}}
   \put(4,0){\line(1,0){12}}
   \put(16,0){\line(0,1){16}}
   \put(8,0){\line(0,1){8}}
   \put(1,-2){7}
   \put(3,7){$B$}
   \put(12,15){8}
  \end{picture}
    }
 \put(165,90)%
    { \begin{picture}(10,10)
   \put(0,0){\circle*{2.7}}
   \put(8,8){\circle*{2}}
   \put(0,0){\line(1,0){8}}
   \put(8,0){\line(0,1){8}}
   \put(4,7){6}
   \put(-4,-2){$b$}
  \end{picture}
    }
 \put(150,110)%
    { \begin{picture}(25,25)
   \put(8,8){\circle*{2.7}}
   \put(16,16){\circle*{2}}
   \put(24,24){\circle*{2}}
   \put(16,16){\line(1,0){8}}
   \put(8,8){\line(1,0){16}}
   \put(4,0){\line(1,0){20}}
   \put(24,0){\line(0,1){24}}
   \put(16,0){\line(0,1){16}}
   \put(8,0){\line(0,1){8}}
   \put(1,-2){3}
   \put(12,15){4}
   \put(3,8){$A$}
   \put(20,23){5}
   \put(0,25){(d)}

  \end{picture}
    }
 \put(135,75)%
    { \begin{picture}(20,20)
   \put(0,0){\circle*{2.7}}
   \put(16,16){\circle*{2}}
   \put(8,8){\circle*{2}}
   \put(8,8){\line(1,0){8}}
   \put(0,0){\line(1,0){16}}
   \put(16,0){\line(0,1){16}}
   \put(8,0){\line(0,1){8}}
   \put(-4,-3){$a$}
   \put(5,7){1}
   \put(12,15){2}
  \end{picture}
    }
 \put(140,125)%
    { \begin{picture}(7,7)
   \put(0,0){\circle*{2.7}}
   \put(2,0){$c$}
  \end{picture}
    }
   \end{picture}
  \end{center}
\caption{Creation of $S(2,1)$. (a) $\Gamma(2,1)$; (b) $C(2,1)$; (c) the
left $B$-sails; (d) the right $B$-sails; (e) $S(2,1)$. }
 \label{fig:fig4b}
  \end{figure}


\begin{figure}[htb]
 \begin{center}
 \unitlength=0.8mm
\begin{picture}(180,175)
 \put(70,140)%
    {\begin{picture}(40,35)
   \put(10,5){\circle*{2.7}}
   \put(25,5){\circle*{2.7}}
   \put(10,20){\circle*{2.7}}
   \put(25,20){\circle*{2.7}}
   \put(10,35){\circle*{2.7}}
   \put(25,35){\circle*{2.7}}
   \put(10,5){\line(1,0){15}}
   \put(10,20){\line(1,0){15}}
   \put(10,35){\line(1,0){15}}
   \put(10,5){\line(0,1){30}}
   \put(25,5){\line(0,1){30}}
   \put(6,2){$a$}
   \put(26.5,2.5){$A$}
   \put(6,18){$b$}
   \put(27,18){$B$}
   \put(6,35){$c$}
   \put(27,35){$C$}
   \put(-5,20){(a)}
     \end{picture}
    }
 \put(60,85)%
     {\begin{picture}(50,45)
   \put(10,5){\circle*{2.7}}
   \put(18,15){\circle*{2.7}}
   \put(25,20){\circle*{2.7}}
   \put(33,30){\circle*{2.7}}
   \put(40,35){\circle*{2.7}}
   \put(48,45){\circle*{2.7}}
   \put(18,5){\circle*{2}}
   \put(10,20){\circle*{2}}
   \put(33,20){\circle*{2}}
   \put(18,30){\circle*{2}}
   \put(10,35){\circle*{2}}
   \put(25,35){\circle*{2}}
   \put(48,35){\circle*{2}}
   \put(18,45){\circle*{2}}
   \put(33,45){\circle*{2}}
   \put(10,5){\line(1,0){8}}
   \put(10,20){\line(1,0){23}}
   \put(18,30){\line(1,0){15}}
   \put(10,35){\line(1,0){38}}
   \put(18,45){\line(1,0){30}}
   \put(10,5){\line(0,1){30}}
   \put(18,5){\line(0,1){40}}
   \put(25,20){\line(0,1){15}}
   \put(33,20){\line(0,1){25}}
   \put(48,35){\line(0,1){10}}
   \put(6,2){$a$}
   \put(26,15){$b$}
   \put(41,31.5){$c$}
   \put(13,13){$A$}
   \put(35,27){$B$}
   \put(50,45){$C$}
   \put(40,10){(b)}
   \put(6,18){1}
   \put(6,34){2}
   \put(20,2){3}
   \put(14,28){4}
   \put(14,44){5}
   \put(25,37){6}
   \put(35,17){7}
   \put(34,46){8}
   \put(50,32){9}
    \end{picture}
    }
 \put(50,0)%
    {\begin{picture}(70,72)
   \put(0,0){\circle*{2.7}}
   \put(30,10){\circle*{2.7}}
   \put(24,31){\circle*{2.7}}
   \put(40,41){\circle*{2.7}}
   \put(34,62){\circle*{2.7}}
   \put(64,72){\circle*{2.7}}
   \put(20,5){\circle*{2}}
   \put(10,50){\circle*{2}}
   \put(54,22){\circle*{2}}
   \put(32,36){\circle*{2}}
   \put(44,67){\circle*{2}}
   \put(20,0){\circle*{1.2}}
   \put(20,10){\circle*{1.2}}
   \put(20,16){\circle*{1.2}}
   \put(30,16){\circle*{1.2}}
   \put(40,16){\circle*{1.2}}
   \put(20,22){\circle*{1.2}}
   \put(30,22){\circle*{1.2}}
   \put(0,26){\circle*{1.2}}
   \put(24,26){\circle*{1.2}}
   \put(32,26){\circle*{1.2}}
   \put(32,31){\circle*{1.2}}
   \put(32,41){\circle*{1.2}}
   \put(32,46){\circle*{1.2}}
   \put(40,46){\circle*{1.2}}
   \put(64,46){\circle*{1.2}}
   \put(0,50){\circle*{1.2}}
   \put(34,50){\circle*{1.2}}
   \put(44,50){\circle*{1.2}}
   \put(24,56){\circle*{1.2}}
   \put(34,56){\circle*{1.2}}
   \put(44,56){\circle*{1.2}}
   \put(44,62){\circle*{1.2}}
   \put(44,72){\circle*{1.2}}
   \put(0,0){\line(1,0){20}}
   \put(20,10){\line(1,0){10}}
   \put(20,16){\line(1,0){20}}
   \put(20,22){\line(1,0){44}}
   \put(0,26){\line(1,0){32}}
   \put(24,31){\line(1,0){8}}
   \put(32,41){\line(1,0){8}}
   \put(32,46){\line(1,0){32}}
   \put(0,50){\line(1,0){44}}
   \put(24,56){\line(1,0){20}}
   \put(34,62){\line(1,0){10}}
   \put(44,72){\line(1,0){20}}
   \put(0,0){\line(0,1){50}}
   \put(20,0){\line(0,1){22}}
   \put(24,26){\line(0,1){30}}
   \put(30,10){\line(0,1){12}}
   \put(32,26){\line(0,1){20}}
   \put(34,50){\line(0,1){12}}
   \put(40,16){\line(0,1){30}}
   \put(44,50){\line(0,1){22}}
   \put(64,22){\line(0,1){50}}
 %
   \put(30,16.3){\line(1,0){10}}
   \put(30,15.7){\line(1,0){10}}
   \put(0,26.3){\line(1,0){24}}
   \put(0,25.7){\line(1,0){24}}
   \put(40,46.3){\line(1,0){24}}
   \put(40,45.7){\line(1,0){24}}
   \put(24,56.3){\line(1,0){10}}
   \put(24,55.7){\line(1,0){10}}
   \put(-4,-3){$a$}
   \put(31,7){$A$}
   \put(20.5,31){$b$}
   \put(41.5,38){$B$}
   \put(30,62){$c$}
   \put(66,72){$C$}
   \put(22,3){3}
   \put(9,52){2}
   \put(53,17){5}
   \put(34,34){7}
   \put(46,65){9}
   \put(10,27.5){$e_1$}
   \put(35,12.5){$e_4$}
   \put(27,57.5){$e_6$}
   \put(50,42.5){$e_8$}
   \put(-15,35){(e)}
  \end{picture}
    }
%
 \put(25,90)%
    { \begin{picture}(10,10)
   \put(0,0){\circle*{2.7}}
   \put(0,8){\line(1,0){4}}
   \put(0,0){\line(0,1){8}}
   \put(2,-3){$b$}
   \put(5,6.5){6}
  \end{picture}
    }
 \put(20,125)%
    { \begin{picture}(20,20)
   \put(8,8){\circle*{2.7}}
   \put(0,0){\circle*{2}}
   \put(0,8){\line(1,0){8}}
   \put(0,16){\line(1,0){12}}
   \put(0,0){\line(0,1){16}}
   \put(8,8){\line(0,1){8}}
   \put(2,-3){7}
   \put(9,5){$B$}
   \put(13,14){8}
  \end{picture}
    }
 \put(30,110)%
    { \begin{picture}(10,10)
   \put(8,8){\circle*{2.7}}
   \put(0,0){\circle*{2}}
   \put(0,8){\line(1,0){8}}
   \put(0,0){\line(0,1){8}}
   \put(2,-3){9}
   \put(9.5,7){$C$}
  \end{picture}
    }
 \put(0,90)%
    { \begin{picture}(25,25)
   \put(8,8){\circle*{2.7}}
   \put(0,0){\circle*{2}}
   \put(16,24){\circle*{2}}
   \put(0,8){\line(1,0){8}}
   \put(0,16){\line(1,0){12}}
   \put(0,24){\line(1,0){16}}
   \put(0,0){\line(0,1){24}}
   \put(8,8){\line(0,1){16}}
   \put(2,-2){3}
   \put(9,5){$A$}
   \put(13,14){4}
   \put(18,23){5}
   \put(2,40){(c)}

  \end{picture}
    }
 \put(5,65)%
    { \begin{picture}(20,20)
   \put(0,0){\circle*{2.7}}
   \put(8,16){\circle*{2}}
   \put(0,8){\line(1,0){4}}
   \put(0,16){\line(1,0){8}}
   \put(0,0){\line(0,1){16}}
   \put(2,-3){$a$}
   \put(4.5,6){1}
   \put(10,15){2}
  \end{picture}
    }
 \put(30,75)%
    { \begin{picture}(7,7)
   \put(0,0){\circle*{2.7}}
   \put(2,0){$c$}
  \end{picture}
    }
%
%
 \put(150,90)%
    { \begin{picture}(10,10)
   \put(8,8){\circle*{2.7}}
   \put(4,0){\line(1,0){4}}
   \put(8,0){\line(0,1){8}}
   \put(3.5,8){$B$}
   \put(1,-2){4}
  \end{picture}
    }
 \put(165,85)%
    { \begin{picture}(20,20)
   \put(8,8){\circle*{2.7}}
   \put(16,16){\circle*{2}}
   \put(8,8){\line(1,0){8}}
   \put(4,0){\line(1,0){12}}
   \put(16,0){\line(0,1){16}}
   \put(8,0){\line(0,1){8}}
   \put(1,-2){1}
   \put(3,7){$b$}
   \put(12.5,15){7}
  \end{picture}
    }
 \put(150,70)%
    { \begin{picture}(10,10)
   \put(0,0){\circle*{2.7}}
   \put(8,8){\circle*{2}}
   \put(0,0){\line(1,0){8}}
   \put(8,0){\line(0,1){8}}
   \put(4,7){3}
   \put(-4,-2){$a$}
  \end{picture}
    }
 \put(140,110)%
    { \begin{picture}(25,25)
   \put(8,16){\circle*{2.7}}
   \put(0,0){\circle*{2}}
   \put(16,24){\circle*{2}}
   \put(0,0){\line(1,0){16}}
   \put(4,8){\line(1,0){12}}
   \put(8,16){\line(1,0){8}}
   \put(8,0){\line(0,1){16}}
   \put(16,0){\line(0,1){24}}
   \put(-3.5,-3){2}
   \put(4.5,16){$c$}
   \put(2,6){6}
   \put(12,23){9}
   \put(0,30){(d)}

  \end{picture}
    }
 \put(170,120)%
    { \begin{picture}(20,20)
   \put(8,16){\circle*{2.7}}
   \put(0,0){\circle*{2}}
   \put(0,0){\line(1,0){8}}
   \put(4,8){\line(1,0){4}}
   \put(8,0){\line(0,1){16}}
   \put(-4,-3){5}
   \put(1,6){8}
   \put(3,15){$C$}
  \end{picture}
    }
 \put(140,85)%
    { \begin{picture}(7,7)
   \put(0,0){\circle*{2.7}}
   \put(2,0){$A$}
  \end{picture}
    }
   \end{picture}
  \end{center}
\caption{Creation of $S(1,2)$. (a) $\Gamma(1,2)$; (b) $C(1,2)$; (c) the
left $B$-sails; (d) the right $B$-sails; (e) $S(1,2)$. }
 \label{fig:fig4c}
  \end{figure}

We will prove later, via a chain of equivalence relations, that the
S-graphs constructed above are precisely the regular $B_2$-crystals.

\medskip
 \noindent
{\bf Remark 1.} We have seen that the principal vertices (i.e.,
those coming from the corresponding grid $\Gamma$) are important
in the construction of $S(a,b)$. They also possess the following
nice property: for any two principal vertices $v,v'$ whose
coordinates $(i,j)$ and $(i',j')$ (respectively) in $\Gamma(a,b)$
satisfy $i\le i'$ and $j\le j'$, the interval of $S(a,b)$ from $v$
to $v'$ is isomorphic to $S(i'-i,j'-j)$. (In an acyclic digraph
$G$, an interval from a vertex $x$ to a vertex $y$ is the subgraph
of $G$ whose vertices and edges belong to paths from $x$ to $y$.)
This property can be deduced from the above construction and
becomes quite transparent when $S(a,b)$ is represented via the
worm model introduced in Section~\ref{sec:worm}. (A similar
property for principal vertices in regular $A_n$-crystals is
shown in~\cite{A-N}. See also Lemma~7.12 in~\cite{Kbook}.)

Another feature that can be obtained from the construction (and will be
easily seen from the worm model; cf. Remark~4 in Section~\ref{sec:worm})
is the existence of a
mapping from the vertex set of $S(a,b)$ to $\Zset^2$ which brings each
1-edge to a vector congruent to $(1,0)$, and each 2-edge to a vector
congruent to $(0,1)$ (a ``weight mapping'').

\medskip
{\bf Infinite S-graphs}.

\smallskip
 In~\cite{A-2} the construction of regular $A_2$-crystals is
generalized in a natural way to produce their infinite analogs. To do
so, one considers a grid $\Gamma=P\times Q$ in which one of the paths
$P,Q$ or both is allowed to be semi-infinite in forward direction or
semi-infinite in backward direction or fully infinite, i.e. to be of the
form $\ldots,v_i,(v_i,v_{i+1}),v_{i+1},(v_{i+1},v_{i+2}),\ldots$ with
the indices running over $\Zset_+$ or $\Zset_-$ or $\Zset$,
respectively. The definitions of right $A$-sails $R$ and left $A$-sails
$L$ are extended accordingly. For example, if $P$ is fully infinite,
then the vertex set of $R$ is formed by the points $(i,j)\in \Zset^2$
with $j\le i$, and if $Q$ is semi-infinite in backward direction, then
the vertex set of $L$ is formed by $(i,j)\in \Zset^2$ with $i\le j\le
0$. As before, the ``diagonals'' in these sails are the sets of pairs
$(i,i)$. The diagonal-product $R$\sprod$L$ (which remains well-defined
in this general case) is the desired graph; following terminology
in~\cite{A-2}, we call such a (finite or infinite) graph $C$ an {\em
RA2-graph}. Now the above construction of right and left $B$-sails is
extended in a natural way to treat I- and II-strings of such a $C$
(note that, as before, the diagonal side of each $B$-sail contains
exactly one (finite) break point). As a result, we obtain an extended
collection which includes the above S-graphs and their infinite
analogs, to which we refer as S-graphs as well. We may conditionally
use symbols $\Zset_+,\Zset_-,\Zset$ to denote corresponding parameters
(in place of finite numbers $a,b$); e.g. $S(a,\Zset_-)$ denotes the
S-graph determined by a finite path $P_a$ with color $\alpha$ and the
semi-infinite in backward direction path with color $\beta$.


\section{Axioms}  \label{sec:axiom}

In this section we present a list of axioms defining a class of
2-edge-colored digraphs. Although finite graphs are of most interest
for us, our axioms are stated so that they fit infinite graphs as well.
We will prove that the graphs defined by these axioms are
exactly the S-graphs from the previous section (thus obtaining an
axiomatical characterization of the regular $B_2$-crystals, by
attracting further arguments).

Let $G=(V,E)$ be a (weakly) connected 2-edge-colored digraph
without multiple edges ($G$ is allowed to be infinite). As before, we
denote the edge colors of $G$ by 1 and 2 and refer to an
edge with color $i$ as an $i$-{\em edge}. 

Our original list consists of Axioms (B0)--(B4), Axioms (B3$'$)--(B4$'$)
(``dual'' of (B3)--(B4)), and Axiom~(BA). The last axiom bridges $B_2$-
and $A_2$-crystals; it requires that a certain graph derived from $G$
(satisfying Axioms (B0)--(B4),(B3$'$)--(B4$'$)) be an RA2-graph in which a
certain partition of the vertices into two subsets is distinguished. Such a ``decorated'' RA2-graph will be characterized via ``local'' axioms
(Axioms~(A0)--(A9) in this section), and we will translate these
axioms in direct terms of the original graph $G$ in Appendix~2.

\begin{itemize}
\item[{\bf (B0)}]
For $i=1,2$, each vertex has at most one entering edge
with color $i$ and at most one leaving edge with color $i$.

 \end{itemize}

Thus, the deletion of all edges of $G$ of color $3-i$ produces a
disjoint union of (finite or infinite) paths, called strings of color
$i$, or $i$-{\em strings}.
The 1-edges (2-edges) can be identified with the action of the
corresponding partial invertible operator on the vertices, denoted
by {\bf 1} (resp. by {\bf 2}); in particular, for a 1-edge $(u,v)$, we
may write $v={\bf 1}u$ and $u={\bf 1}^{-1}v$.

In the next two axioms we assume that some vertices and 1-edges of $G$ are distinguished as \emph{central} ones. 

\begin{itemize}
\item[{\bf (B1)}]
Each 1-string has exactly one central element, which is either a vertex or an edge.
  \end{itemize}

Due to this axiom, the vertex set $V$ of $G$ is partitioned into three
subsets: the sets of central, left and right vertices. Here a vertex
$v$ is called {\em left} ({\em right}) if it lies before (resp. after)
the central element of the 1-string containing $v$. Accordingly, a
non-central 1-edge $(u,v)$ is regarded as left (right) if $u$ is left
(resp. $v$ is right). Note that for a central edge $e=(x,y)$, the
vertex $x$ is left and the vertex $y$ is right; it will be convenient
for us to think of $e$ as though consisting of two {\em half-edges}:
the left half-edge, beginning at $x$ and ending at the ``mid-point'' of
$e$, and the right half-edge, beginning at the ``mid-point'' and ending
at $y$.

\begin{itemize}
\item[{\bf (B2)}]
Each 2-string $P$ contains exactly one central vertex $v$. Moreover,
all vertices of $P$ lying {\em before} $v$ are right, whereas all
vertices lying {\em after} $v$ are left.
  \end{itemize}

\noindent(So there is a difference in the locations of left and right vertices in 1- and 2-strings.) A 2-edge $(u,v)$ is regarded as left (right) if $v$ is left (resp. $u$
is right).

\begin{corollary}  \label{cor:enter-young}
For any left or central 1-edge $(u,v)$, there exists a 2-edge entering
$u$, and the latter edge is left. For any left 2-edge $(u',v')$, there
exists a 1-edge leaving $v'$, and at least a half of the latter edge is
left.
  \end{corollary}

Indeed, consider the 2-string $P$ passing $u$. Since $u$ is left, the
center of $P$ is located strictly before $u$, by (B2). This gives the
first assertion, and the second one is shown in a similar way.

Strictly speaking, Axioms (B1) and (B2) are not local. In order to obtain their local implementations, we assume that, instead of distingushing central elements as before, the graph $G=(V,E)$ is equipped with \emph{labels} $\ell(v)$ on the vertices $v\in V$ which take values in the 3-element set $\{0,c,1\}$. Assuming that $G$ is acyclic (i.e. has no directed cycle), we impose the following local requirements on $(G,\ell)$.

 \begin{itemize}
  \item[(B1(i))]  For each 1-edge $(u,v)$, the pair $(\ell(u),\ell(v))$ is equal to one of $(0,0),(0,c),(0,1)$, $(c,1),(1,1)$.
  
  \item[(B1(ii))] If a vertex $v$ has no entering 1-edge, then $\ell(v)\ne 1$ (i.e., $\ell(v)\in\{0,c\}$); and if a vertex $v$ has no leaving 1-edge, then $\ell(v)\ne 0$.
  
  \item[(B2(i))] For each 2-edge $(u,v)$, the pair $(\ell(u),\ell(v))$ is equal to one of $(1,1),(1,c)$, $(c,0),(0,0)$.  
  
  \item[(B2(ii))] If a vertex $v$ has no entering 2-edge, then $\ell(v)\ne 0$; and if a vertex $v$ has no leaving 2-edge, then $\ell(v)\ne 1$.  
  \end{itemize}

\noindent (In particular, (B1(ii)) implies that if $v$ has neither entering nor leaving 1-edge, then $\ell(v)=c$, and similarly for (B2(ii)).)

The vertices labeled $0,c,1$ are naturally interpreted as left, central and right ones, respectively. Accordingly, a 1-edge $(u,v)$ is regarded as left if $(\ell(u),\ell(v))\in\{(0,0),(0,c)\}$, central if $(\ell(u),\ell(v))=(0,1)$, and right if $(\ell(u),\ell(v))\in\{(c,1),(1,1)\}$. As to a 2-edge $(u,v)$, it is regarded as left if $(\ell(u),\ell(v))\in\{(c,0),(0,0)\}$, and right if $(\ell(u),\ell(v))\in\{(1,1),(1,c)\}$.  

As an easy consequence of the above assignments, the following takes place:

\begin{itemize} 
\item \emph{For finite acyclic graphs satisfying~(B0), Axioms (B1)--(B2) are equivalent to imposing (B1(i),(ii)),(B2(i),(ii)).}
 \end{itemize}

Based on this equivalence, it will eventually follow that when dealing with finite acyclic graphs, we will obtain a ``purely local'' axiomatics for finite $B_2$-crystals. Note also that the formal requirement that an input graph $G$ is acyclic can be replaced by imposing additional local variables and constraints. Namely, we can endow each vertex $v$ with an additional number $\pi(v)\in\Zset$ (a ``potential'') and impose the condition: $\pi(u)<\pi(v)$ for each edge $(u,v)$ of $G$. Then $G$ is acyclic if and only if a feasible $\pi$ exists.

However, the above labeling method does not work when an input graph is infinite (since in this case Axioms (B1(i),(ii)),(B2(i),(ii)) do not forbid infinite monochromatic strings without central elements). Therefore, our ``local implementation'' is applicable only to finite $B_2$-crystals.
\medskip

For Axioms (B3),(B4) stated below (as well as for Axioms (B5)--(B13)
given in Appendix~2), the corresponding dual axiom will also be default
imposed on $G$; this is the same statement applied to the {\em dual
graph}, the graph obtained from $G$ by reversing the orientation of all
edges while preserving their colors and distinguishing the same central
elements (so, to get a formulation for $G$ itself, one should swap the
words ``left'' and ``right'', as well as ``enter'' and ``leave''). The
axiom dual of (Bi) is denoted by (Bi$'$).

In illustrations below we indicate the central vertices by thick dots
and mark the central edges by black rhombi in the middle of the
corresponding arrows.

By a {\em commutative square} one means a quadruple of vertices
$u,v,u',v'$ (in any order) such that $v={\bf 1}u$, $u'={\bf 2}u$ and
$v'={\bf 1}u'={\bf 2}v$, as well as the subgraph of $G$ induced by
these vertices. The next axiom (together with its dual) says that a
non-commutativity of the operators {\bf 1} and {\bf 2} may occur only
in a vicinity of central vertices or central edges.

\begin{itemize}
\item[{\bf (B3)}]
(i) Let $(u,v)$ be a left 1-edge and suppose that there is a 2-edge
leaving $u$ or $v$. Then these two edges belong to a commutative
square. (ii) Let $(v,w)$ be a left 2-edge and suppose that there is a
1-edge entering $v$ or $w$. Then these two edges belong to a
commutative square. (See the picture.)
  \end{itemize}

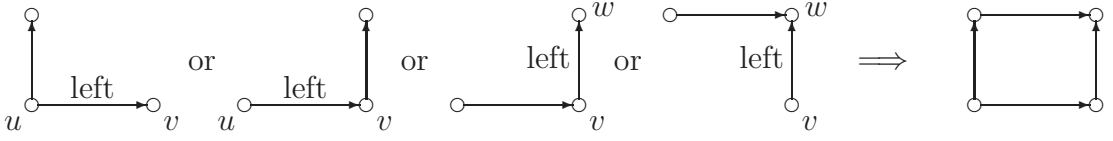
\unitlength=.8mm \special{em:linewidth 0.4pt} \linethickness{0.4pt}
\begin{picture}(185.00,25.00)(0,0)
  \put(5.00,5.00){\circle{2.00}}
 \put(25.00,5.00){\circle{2.00}}
 \put(5.00,20.00){\circle{2.00}}
 \put(6.00,5.00){\vector(1,0){18.00}}
 \put(5.00,6.00){\vector(0,1){13.00}}
 \put(2.00,2.00){\makebox(0,0)[cc]{$u$}}
 \put(28.00,2.00){\makebox(0,0)[cc]{$v$}}
 \put(15.00,8.00){\makebox(0,0)[cc]{left}}
\put(33,12){\makebox(0,0)[cc]{or}}
  \put(40.00,5.00){\circle{2.00}}
  \put(60.00,5.00){\circle{2.00}}
 \put(60.00,20.00){\circle{2.00}}
 \put(41.00,5.00){\vector(1,0){18.00}}
 \put(60.00,6.00){\vector(0,1){13.00}}
 \put(37.00,2.00){\makebox(0,0)[cc]{$u$}}
 \put(63.00,2.00){\makebox(0,0)[cc]{$v$}}
 \put(50.00,8.00){\makebox(0,0)[cc]{left}}
 \put(68,12){\makebox(0,0)[cc]{or}}
  \put(75.00,5.00){\circle{2.00}}
  \put(95.00,5.00){\circle{2.00}}
 \put(95.00,20.00){\circle{2.00}}
 \put(76.00,5.00){\vector(1,0){18.00}}
 \put(95.00,6.00){\vector(0,1){13.00}}
 \put(98.00,2.00){\makebox(0,0)[cc]{$v$}}
 \put(99.00,21.00){\makebox(0,0)[cc]{$w$}}
 \put(90.00,13.00){\makebox(0,0)[cc]{left}}
\put(103,12){\makebox(0,0)[cc]{or}}
  \put(110.00,20.00){\circle{2.00}}
  \put(130.00,5.00){\circle{2.00}}
 \put(130.00,20.00){\circle{2.00}}
 \put(111.00,20.00){\vector(1,0){18.00}}
 \put(130.00,6.00){\vector(0,1){13.00}}
 \put(133.00,2.00){\makebox(0,0)[cc]{$v$}}
 \put(134.00,21.00){\makebox(0,0)[cc]{$w$}}
 \put(125.00,13.00){\makebox(0,0)[cc]{left}}
\put(145,12){\makebox(0,0)[cc]{$\Longrightarrow$}}
 \put(160.00,5.00){\circle{2.00}}
 \put(180.00,5.00){\circle{2.00}}
 \put(180.00,20.00){\circle{2.00}}
 \put(160.00,20.00){\circle{2.00}}
 \put(161.00,5.00){\vector(1,0){18.00}}
 \put(160.00,6.00){\vector(0,1){13.00}}
 \put(161.00,20.00){\vector(1,0){18.00}}
 \put(180.00,6.00){\vector(0,1){13.00}}
  \end{picture}

(Then the first part of the dual axiom (B3$'$) says that for a right 1-edge
$(u,v)$, if $u$ or $v$ has entering 2-edge, then these two edges belong
to a commutative square.)

Note that the square $uu'vv'$ as in (B3)(i), where $u'={\bf 2}u$ and
$v'={\bf 2}v$, is commutative in a stronger sense:
\smallskip 

$(\ast)$ ~\emph{the edges $(u,u')$,
$(v,v')$, $(u',v')$ are left as well}. 
\smallskip

Indeed, the vertex $v$ is left or
central. Therefore, the vertex $v'$ is left, by (B2), whence the above
edges are left.
 \medskip

 \noindent
{\bf Remark 2.} We shall see later that when $G$ is finite, the
difference between central, left and right edges consists in the
following. For a vertex $v$ and color $i$, let $t_i(v)$ (resp.
$h_i(v)$) denote the length of the part of the $i$-string containing
$v$ from the beginning to $v$ (resp. from $v$ to the end). For an
$i$-edge $e=(u,v)$, define $\Delta t(e):=t_{3-i}(v)-t_{3-i}(u)$ and
$\Delta h(e):=h_{3-i}(v)-h_{3-i}(u)$. Then the status of $e$ is
described in terms of the pair $\Delta(e):=(\Delta t(e),\Delta h(e))$:
(i) when $i=1$, $\Delta(e)$ is $(-2,0)$ if $e$ is left, $(0,2)$ if $e$
is right, and $(-1,1)$ if $e$ is central; (ii) when $i=2$, $\Delta(e)$
is $(-1,0)$ if $e$ is right, and $(0,1)$ if $e$ is left. (Note that the
value $\Delta t(e)-\Delta h(e)$ is equal to $-2$ for all 1-edges $e$
and $-1$ for all 2-edges $e$, which matches the off-diagonal
coefficients of the Cartan matrix for $B_2$.)
 \medskip
 
The next axiom shows a relation between central edges.

\begin{itemize}
\item[{\bf (B4)}]
Let $u$ be the beginning vertex of a central edge. Suppose that there are a
2-edge $(u,w)$, a 1-edge $(w,w')$ and a 2-edge $(w',v)$. Then $v$ is
the beginning vertex of a central edge.
  \end{itemize}

See the picture (where in the right fragment, the vertex $w'$ is indicated as central due to Corollary~\ref{cor:centr-edge2}).

\unitlength=.8mm \special{em:linewidth 0.4pt} \linethickness{0.4pt}
\begin{picture}(120.00,40.00)(0,0)
  {
  \begin{picture}(71.00,40.00)(10,0)
 \put(40.00,5.00){\circle{2.00}}
 \put(40.00,20.00){\circle{2.00}}
 \put(60.00,20.00){\circle{2.00}}
 \put(60.00,35.00){\circle{2.00}}
 \put(49.00,5.00){\makebox(0,0)[cc]{\scriptsize $\blacklozenge$}}
 \put(40.00,6.00){\vector(0,1){13.00}}
 \put(41.00,20.00){\vector(1,0){18.00}}
 \put(60.00,21.00){\vector(0,1){13.00}}
 \put(41.00,5.00){\vector(1,0){9.00}}
 \put(64.00,18.00){\makebox(0,0)[cc]{$w'$}}
 \put(36.00,19.00){\makebox(0,0)[cc]{$w$}}
 \put(36.00,3.00){\makebox(0,0)[cc]{$u$}}
 \put(57.00,37.00){\makebox(0,0)[cc]{$v$}}
  \end{picture}
\put(5,18){\makebox(0,0)[cc]{$\Longrightarrow$}}
  \begin{picture}(71.00,40.00)(10,0)
\put(40.00,5.00){\circle{2.00}} \put(40.00,20.00){\circle{2.00}}
\put(60.00,20.00){\circle*{2.00}} \put(60.00,35.00){\circle{2.00}}
\put(69.00,35.00){\makebox(0,0)[cc]{\scriptsize $\blacklozenge$}}
\put(49.00,5.00){\makebox(0,0)[cc]{\scriptsize $\blacklozenge$}}

\put(40.00,6.00){\vector(0,1){13.00}}
\put(41.00,20.00){\vector(1,0){18.00}}
\put(60.00,21.00){\vector(0,1){13.00}}
\put(61.00,35.00){\vector(1,0){9.00}}
\put(41.00,5.00){\vector(1,0){9.00}}
\put(64.00,18.00){\makebox(0,0)[cc]{$w'$}}
\put(36.00,19.00){\makebox(0,0)[cc]{$w$}}
\put(36.00,3.00){\makebox(0,0)[cc]{$u$}}
\put(57.00,37.00){\makebox(0,0)[cc]{$v$}}
  \end{picture}
  }
\end{picture}

\begin{corollary}  \label{cor:centr-edge}
Let $(u,v)$ be a central 1-edge. Then there are a 2-edge $(u',u)$ and a
2-edge $(v,v')$. Moreover, both vertices $u',v'$ are central.
  \end{corollary}

 \begin{proof}
~Since the edge $(u,v)$ is central, the vertex $u$ is left and the vertex $v$ is right. Hence
 $u$ has an entering 2-edge, $(u',u)$ say, and $v$ has a leaving 2-edge, $(v,v')$ say; note that the former edge is left and the latter edge is right. Suppose that the vertex $u'$ is not central. Then $u'$ is left, and therefore $u'$ has a leaving 1-edge, $(u',u'')$ say. The edge $(u',u'')$ cannot be left; for otherwise $u',u'',u,v$ would give a commutative square (by Axiom (B3)), whence $(u,v)$ should be left, not central. Therefore, $(u',u'')$ must be central. But then, applying Axiom (B4) to the edges $(u',u''),(u',u),(u,v),(v,v')$, we observe that the vertex $v'$ is left; this implies that the vertex $v$ is left or central, not right, contradicting the condition that the edge $(u,v)$ is central. Thus, $u'$ is central, as required.
 
The assertion that $v'$ is central is symmetric.
  \end{proof}

\begin{corollary}  \label{cor:centr-edge2}
Let $(u,v)$ be a central 1-edge. Let $u$ have a leaving 2-edge $(u,w)$. Then
there is a 1-edge $(w,w')$, and the vertex $w'$ is central.
  \end{corollary}

 \begin{proof}
~Since the vertex $u$ is left, $w$ is left as well and a 1-edge $(w,w')$ does exist. This edge is not central; otherwise the vertex $u$ would be central, by Corollary~2, contradicting the condition that $(u,v)$ is central. Suppose that $w'$ is not central. Then $w'$ is left, and therefore, $w'$ has an entering 2-edge, $(w'',w')$ say. This edge is left and we can apply~(B3)(ii) to the edges $(w'',w')$ and $(w,w')$, obtaining a commutative square containing the vertices $w,w',w''$. This square must contain the vertices $u,v$ as well, implying $v=w''$ and leading to a contradiction to the fact that the edge $(u,v)$ is central.
  \end{proof}

Axioms (B0)--(B4),(B3$'$)--(B4$'$) enable us to extract $B$-sails in $G$
(defined in the previous section), as follows. Let us cut $G$ in the
SW--NE direction at each central vertex $v$, that is, split $v$ into
two vertices $v',v''$, making $v'$ incident (in place of $v$) with the
1-edge entering $v$ and the 2-edge leaving $v$, and making $v''$
incident with the 1-edge leaving $v$ and the 2-edge entering $v$, when
such edges exist. Also split each central edge $e=(u,w)$ into two
``half-edges'' $(u,v'_e)$ and $(v''_e,w)$, where $v'_e,v''_e$ are
copies of the mid-point $v_e$ of $e$. By (B0)--(B2), the obtained graph
has two sorts of components: those containing left 1- and 2-edges and
left half-edges, called {\em left components}, and those containing
right 1- and 2-edges and right half-edges, called {\em right
components}. (Note that left or right components consisting of a single
vertex are possible.) Moreover, each monochromatic string of $G$
becomes split (at its central point) into two parts, one lying in a
left component, and the other in a right component. Figure~\ref{fig:fig5}
illustrates the splitting procedure.

\begin{figure}[htb]
 \begin{center}
\unitlength=0.8mm \special{em:linewidth 0.4pt} \linethickness{0.4pt}
\begin{picture}(150.00,80)(10,10)
 \put(10.00,25.00){\circle{2.00}}
 \put(10.00,43.00){\circle{2.00}}
 \put(10.00,8.00){\vector(0,1){16.00}}
 \put(10.00,26.00){\vector(0,1){17.00}}
 \put(11.00,43.00){\vector(1,0){15.00}}
 \put(11.00,25.00){\vector(1,0){9.00}}
 \put(10.00,7.00){\circle*{2}}
 \put(27.00,43.00){\circle*{2}}
 \put(21.00,25.00){\makebox(0,0)[cc]{\scriptsize $\blacklozenge$}}
 \put(5.00,65.00){\circle{2.00}}
 \put(5.00,48.00){\vector(0,1){16.00}}
 \put(6.00,65.00){\vector(1,0){9.00}}
 \put(5.00,47.00){\circle*{2}}
 \put(16.00,65.00){\makebox(0,0)[cc]{\scriptsize $\blacklozenge$}}
 \put(10.00,83.00){\circle*{2}}
 \put(155.00,65.00){\circle{2.00}}
 \put(155.00,43.00){\circle{2.00}}
 \put(139.00,43.00){\vector(1,0){15.00}}
 \put(155.00,44.00){\vector(0,1){20.00}}
 \put(155.00,66.00){\vector(0,1){16.00}}
 \put(146.00,65.00){\vector(1,0){8.00}}
 \put(138.00,43.00){\circle*{2}}
 \put(155.00,83.00){\circle*{2}}
 \put(146.00,65.00){\makebox(0,0)[cc]{\scriptsize $\blacklozenge$}}
 \put(164.00,23.00){\circle{2.00}}
 \put(164.00,24.00){\vector(0,1){20.00}}
 \put(155.00,23.00){\vector(1,0){8.00}}
 \put(164.00,45.00){\circle*{2}}
 \put(154.00,23.00){\makebox(0,0)[cc]{\scriptsize $\blacklozenge$}}
 \put(155.00,7.00){\circle*{2.5}}
  \begin{picture}(120.00,52.00)(10,-20)
 \put(80,2){\circle*{2}}
 \put(80,14){\circle{2}}
 \put(80,26){\circle{2}}
 \put(95,14){\circle{2}}
 \put(95,38){\circle{2}}
 \put(110,26){\circle{2}}
 \put(110,38){\circle{2}}
 \put(110,50){\circle*{2}}
 \put(91,26){\circle*{2}}
 \put(95,30){\circle*{2}}
\put(81,14){\vector(1,0){13}}
 \put(81,26){\vector(1,0){9}}
 \put(92,26){\vector(1,0){17}}
 \put(96,38){\vector(1,0){13}}
 \put(80,15){\vector(0,1){10}}
 \put(95,15){\vector(0,1){14}}
 \put(95,31){\vector(0,1){6}}
 \put(110,27){\vector(0,1){10}}
 \put(80,3){\vector(0,1){10}}
 \put(110,39){\vector(0,1){10}}
  \put(87,14){\makebox(0,0)[cc]{\scriptsize $\blacklozenge$}}
  \put(102,38){\makebox(0,0)[cc]{\scriptsize $\blacklozenge$}}
    \end{picture}
   \end{picture}
  \end{center}
\caption{Splitting the graph $S(0,2)$ (in the middle) into left and
right components.}
\label{fig:fig5}
  \end{figure}

Consider a left component $K$. By the commutativity axiom (B3),
``above'' any 1-edge (not a half-edge) and ``to the left'' of any
2-edge of $K$ we have a strip of commutative squares of $G$; they
belong to $K$ as well. Therefore, $K$ forms a region ``without holes''
of the grid on $\Zset^2$ and has 3 boundaries: the left boundary,
possibly at infinity, formed by the beginnings of the 1-strings
intersecting $K$, the upper boundary, possibly at infinity, formed by
the ends of the 2-strings intersecting $K$, and the lower-right
boundary formed by the central points (central vertices and mid-points
of central edges) occurring in $K$; we call the third boundary, $D(K)$,
the {\em diagonal} of $K$. When the left (upper) boundary is not at
infinity, it is the left part of a 2-string (resp. the left part of a
1-string). To show that $K$ is a (possibly infinite) left $B$-sail, we
have to examine the diagonal $D(K)$.

We order the elements of $D(K)$ in a natural way, by increasing the
coordinates in the grid.  Axiom~(B3) and
Corollary~\ref{cor:enter-young} imply that for each central vertex
$v\in D(K)$, the next element is either the central vertex of the form
${\bf 1}{\bf 2}v$ or the mid-point of the central edge $({\bf 2}v,{\bf
1}{\bf 2}v)$. Let $D(K)$ contain the mid-point $v'_e$ of some central
edge $e=(x,x')$. By Corollary~\ref{cor:centr-edge}, $v'_e$ is not the
first element of $D(K)$, and its preceding element is the central
vertex of the form ${\bf 2}^{-1}x$. Suppose $v'_e$ is not the last
element of $D(K)$. By Corollary~\ref{cor:centr-edge2}, the next element
$y$ is the central vertex of the form ${\bf 1}{\bf 2}x$. Now if $y$ is
not the last in $D(K)$, then, by Axiom~(B4), the next element is again
a mid-point, namely, the mid-point of the central edge $(z,z')$ for
$z={\bf 2}y$.

These observations show that $K$ matches the construction of a
(possibly infinite) left $B$-sail in Section~\ref{sec:constr}, whose
diagonal has a (possibly degenerate) slope 1 followed by a (possibly
degenerate) slope 2. (It is still possible that $D(K)$ is unbounded
from above and has only slope 1, or is unbounded from below and has
only slope 2, but this is not the case as we shall see later.)

Similarly (using the dual statements of
Corollaries~\ref{cor:enter-young},\ref{cor:centr-edge2}), each right
component forms a right $B$-sail. Now we define the following
2-edge-colored digraph $C(G)$, called the {\em central graph} for $G$.
Its vertices are the central vertices and the mid-points of central
edges of $G$. Two vertices $u,v$ of $C(G)$ are connected by an edge,
from $u$ to $v$, of color I if and only if (a copy of) $v$ is next to
$u$ in the diagonal of a right component, and are connected by an edge
of color II if and only if $v$ is next to $u$ in the diagonal of a
left component. (A priori multiple edges are possible.) We also
distinguish between the mid-points, referring to them as
$\otimes$-{\em vertices}, and the other vertices of $C(G)$. Note
that from the connectedness of $G$ it easily follows that $C(G)$ is
connected as well.

Given a (finite or infinite) RA2-graph $C$, let us say that a vertex $v$
of $C$ is an $\otimes$-vertex if $v$ belongs to a {\em left} $A$-sail $LA$
of $C$ and lies at {\em odd} distance from the diagonal of $LA$. The
non-$\otimes$-vertices of $H$ are called {\em ordinary}.

Our final axiom is the following:

\begin{itemize}
\item[{\bf(BA)}] The central graph $C(G)$ is an RA2-graph, and moreover, the
sets of $\otimes$-vertices in these graphs are the same.
  \end{itemize}

\begin{theorem}  \label{tm:ax-S}
Axioms (B0)--(B4),(B3\,$'$)--(B4\,$'$),(BA) define precisely the set of
(finite and infinite) S-graphs.
  \end{theorem}
  \begin{proof}
Immediate from the construction of S-graphs and reasonings above.
  \end{proof}

In the rest of this section we give defining axioms (namely, (A0)--(A8)) for
the RA2-graphs with the distinguished set of $\otimes$-vertices, which we
call the {\em decorated} RA2-graphs. Compared with
``pure'' RA2-graphs, which can be defined via 4 axioms (``local'' or ``almost
local'' ones), cf.~\cite{A-2}, the list of axioms becomes longer because we
wish to describe the difference between $\otimes$-vertices and ordinary ones in
local terms (in the situation when the left sails in an RA2-graph are not
indicated explicitly).

Let $C$ be a connected 2-edge-colored digraph with edge colors I and II
and with a partition of the vertex set into two subsets of which
elements are called $\otimes$-vertices and ordinary vertices,
respectively. In illustrations below the former and latter ones are
indicated by crossed (``otimes'') and white circles, respectively, and as before,
edges of color I are drawn by horizontal arrows.

\begin{itemize}
\item[{\bf(A0)}]
(i) The subgraph of $C$ induced by I-edges consists of pairwise
disjoint paths, and similarly for the II-edges. (ii) No pair of
$\otimes$-vertices is connected by edge. (iii) If $C$ has no
$\otimes$-vertices, then each I-string has a beginning vertex.
(iv) If no pair of ordinary vertices is connected by edge, then
each I-string has an end vertex.
  \end{itemize}

\begin{itemize}
\item[{\bf(A1)}]
Each $\otimes$-vertex has leaving I-edge and entering II-edge.
  \end{itemize}
   \begin{center}
  \unitlength=0.8mm
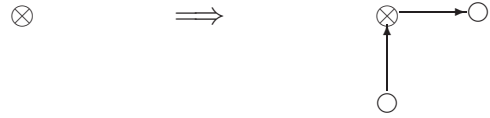
    \begin{picture}(100,22)
 \put(8,18){$\otimes$}
 \put(35,18){$\Longrightarrow$}
 \put(68,18){$\otimes$}
 \put(70,5){\circle{3.0}}
 \put(85,20){\circle{3.0}}
 \put(72,20){\vector(1,0){11}}
 \put(70,7){\vector(0,1){11}}
  \end{picture}
 \end{center}

\begin{itemize}
\item[{\bf(A2)}]
Let an $\otimes$-vertex $v$ have entering I-edge. Then this edge and
the II-edge entering $v$ belong to a square. Moreover, the vertex ${\rm
I}^{-1}{\rm II}^{-1}v$ is an $\otimes$-vertex. (See the picture.)
Symmetrically: if an $\otimes$-vertex $v$ has leaving II-edge, then the
edges leaving $v$ belong to a square, and ${\rm I}({\rm II}v$) is an
$\otimes$-vertex.
  \end{itemize}
   \begin{center}
  \unitlength=0.8mm
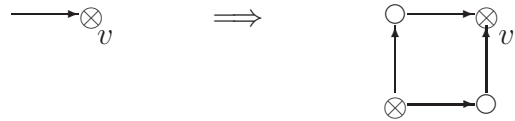
    \begin{picture}(100,22)
 \put(18,18){$\otimes$}
 \put(7,20){\vector(1,0){11}}
 \put(21,15){$v$}
 \put(40,18){$\Longrightarrow$}
 \put(68,3){$\otimes$}
 \put(70,20){\circle{3.0}}
 \put(85,5){\circle{3.0}}
 \put(83,18){$\otimes$}
 \put(72,5){\vector(1,0){11}}
 \put(72,20){\vector(1,0){11}}
 \put(70,7){\vector(0,1){11}}
 \put(85,7){\vector(0,1){11}}
 \put(87,15){$v$}
  \end{picture}
 \end{center}

\begin{itemize}
\item[{\bf(A3)}]
Let an $\otimes$-vertex $v$ have entering I-edge $e=(u,v)$. Suppose
there is a II-edge $e'$ leaving $u$ or $v$. Then $e,e'$ belong to a
square. Moreover, ${\rm II}u$ is an $\otimes$-vertex. (See the
picture.) Symmetrically: if an $\otimes$-vertex $v$ has leaving II-edge
$e=(v,w)$ and if there is a I-edge $e'$ entering $v$ or $w$, then the
edges $e,e'$ belong to a square, and ${\rm I}^{-1}w$ is an
$\otimes$-vertex.
  \end{itemize}
   \begin{center}
  \unitlength=0.8mm
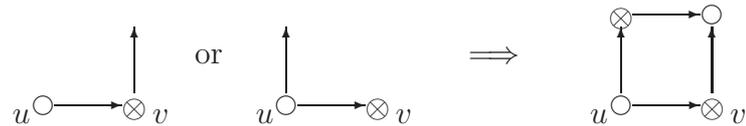
    \begin{picture}(120,22)
 \put(18,3){$\otimes$}
 \put(5,5){\circle{3.0}}
 \put(7,5){\vector(1,0){11}}
 \put(20,7){\vector(0,1){11}}
 \put(0,2){$u$}
 \put(23,2){$v$}
 \put(30,12){or}
 \put(58,3){$\otimes$}
 \put(45,5){\circle{3.0}}
 \put(47,5){\vector(1,0){11}}
 \put(45,7){\vector(0,1){11}}
 \put(40,2){$u$}
 \put(63,2){$v$}
 \put(75,12){$\Longrightarrow$}
 \put(98,18){$\otimes$}
 \put(100,5){\circle{3.0}}
 \put(115,20){\circle{3.0}}
 \put(113,3){$\otimes$}
 \put(102,5){\vector(1,0){11}}
 \put(102,20){\vector(1,0){11}}
 \put(100,7){\vector(0,1){11}}
 \put(115,7){\vector(0,1){11}}
 \put(95,2){$u$}
 \put(118,2){$v$}
  \end{picture}
 \end{center}

\begin{itemize}
\item[{\bf(A4)}]
Let $v$ be an $\otimes$-vertex, and $e=(v,w)$ its leaving I-edge.
Suppose $w$ has leaving II-edge $e'$. Then $e,e'$ belong to a square.
Moreover, ${\rm II}w$ is an $\otimes$-vertex. (See the picture.)
Symmetrically: if $e=(u,v)$ is the II-edge entering an $\otimes$-vertex
$v$ and if $u$ has entering I-edge $e'$, then $e,e'$ belong to a
square, and ${\rm I}^{-1}u$ is an $\otimes$-vertex.
  \end{itemize}
   \begin{center}
  \unitlength=0.8mm
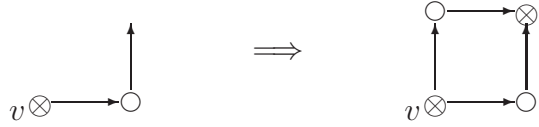
    \begin{picture}(100,22)
 \put(3,3){$\otimes$}
 \put(20,5){\circle{3.0}}
 \put(7,5){\vector(1,0){11}}
 \put(20,7){\vector(0,1){11}}
 \put(0,2){$v$}
 \put(22,2){$w$}
 \put(40,12){$\Longrightarrow$}
 \put(68,3){$\otimes$}
 \put(70,20){\circle{3.0}}
 \put(85,5){\circle{3.0}}
 \put(83,18){$\otimes$}
 \put(72,5){\vector(1,0){11}}
 \put(72,20){\vector(1,0){11}}
 \put(70,7){\vector(0,1){11}}
 \put(85,7){\vector(0,1){11}}
 \put(65,2){$v$}
  \end{picture}
 \end{center}

\begin{itemize}
\item[{\bf(A5)}]
If $(u,v)$ is a I-edge connecting ordinary vertices, then $v$ has
leaving II-edge $(v,w)$, and the vertex $w$ is ordinary. (See the
picture.) Symmetrically: if $(v,w)$ is a II-edge connecting ordinary
vertices, then $v$ has entering I-edge $(u,v)$, and the vertex $u$ is
ordinary.
  \end{itemize}
   \begin{center}
  \unitlength=0.8mm
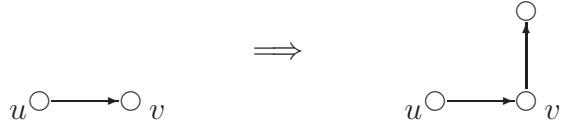
    \begin{picture}(100,22)
 \put(5,5){\circle{3.0}}
 \put(20,5){\circle{3.0}}
 \put(7,5){\vector(1,0){11}}
 \put(0,2){$u$}
 \put(23,2){$v$}
 \put(40,12){$\Longrightarrow$}
 \put(70,5){\circle{3.0}}
 \put(85,5){\circle{3.0}}
 \put(85,20){\circle{3.0}}
 \put(72,5){\vector(1,0){11}}
 \put(85,7){\vector(0,1){11}}
 \put(65,2){$u$}
 \put(88,2){$v$}
  \end{picture}
 \end{center}

\begin{itemize}
\item[{\bf(A6)}]
Let $e$ be a I-edge connecting ordinary vertices $u,v$. Suppose there
is a II-edge $e'$ entering $u$ or $v$. Then $e,e'$ belong to a square,
and all vertices of this square are ordinary. (See the picture.)
Symmetrically: if $e$ is a II-edge connecting ordinary vertices $u,v$
and if $e'$ is a I-edge leaving $u$ or $v$, then $e,e'$ belong to a
square, and all vertices of this square are ordinary.
  \end{itemize}
   \begin{center}
  \unitlength=0.8mm
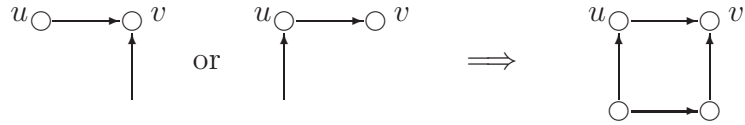
    \begin{picture}(120,22)
 \put(20,20){\circle{3.0}}
 \put(5,20){\circle{3.0}}
 \put(7,20){\vector(1,0){11}}
 \put(20,7){\vector(0,1){11}}
 \put(0,20){$u$}
 \put(23,20){$v$}
 \put(30,12){or}
 \put(60,20){\circle{3.0}}
 \put(45,20){\circle{3.0}}
 \put(47,20){\vector(1,0){11}}
 \put(45,7){\vector(0,1){11}}
 \put(40,20){$u$}
 \put(63,20){$v$}
 \put(75,12){$\Longrightarrow$}
 \put(100,20){\circle{3.0}}
 \put(100,5){\circle{3.0}}
 \put(115,20){\circle{3.0}}
 \put(115,5){\circle{3.0}}
 \put(102,5){\vector(1,0){11}}
 \put(102,20){\vector(1,0){11}}
 \put(100,7){\vector(0,1){11}}
 \put(115,7){\vector(0,1){11}}
 \put(95,20){$u$}
 \put(118,20){$v$}
  \end{picture}
 \end{center}

We say that eight vertices $v_1,\ldots,v_8$ form the {\em small Verma
configuration} from $v_1$ to $v_8$ if, up to renumbering
$v_2,\ldots,v_7$, one has: $v_2={\rm I}v_1$, $v_3={\rm II}v_1$,
$v_4={\rm II}v_2$, $v_5={\rm II}v_4$, $v_6={\rm I}v_3$, $v_7={\rm
I}v_6$, and $v_8={\rm I}v_5={\rm II}v_7$, and in addition,
$v_4\ne v_6$.

\begin{itemize}
\item[{\bf(A7)}]
Let an ordinary vertex $u$ have leaving I-edge $(u,v)$ and leaving
II-edge $(u,w)$. Let $v$ be an ordinary vertex, and $w$ an
$\otimes$-vertex. Then $u,v,w$ belong to the small Verma configuration
from $u$. Moreover, this configuration contains exactly two
$\otimes$-vertices, namely, $w$ and ${\rm II}^2v$.(See the picture.)
Symmetrically: if an $\otimes$-vertex $u$ and ordinary vertices $v,w$
are connected by I-edge $(u,v)$ and II-edge $(w,v)$, then $u,v,w$
belong to the small Verma configuration to $v$. Moreover, this
configuration contains exactly two $\otimes$-vertices, namely, $u$ and
${\rm I}^{-2}w$.
  \end{itemize}
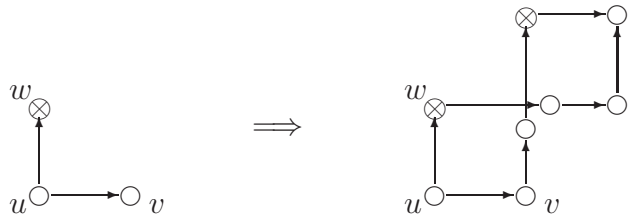
   \begin{center}
  \unitlength=0.8mm
    \begin{picture}(110,40)
 \put(3,18){$\otimes$}
 \put(5,5){\circle{3.0}}
 \put(20,5){\circle{3.0}}
 \put(7,5){\vector(1,0){11}}
 \put(5,7){\vector(0,1){11}}
 \put(0,2){$u$}
 \put(23,2){$v$}
 \put(0,21){$w$}
 \put(40,15){$\Longrightarrow$}
 \put(68,18){$\otimes$}
 \put(70,5){\circle{3.0}}
 \put(85,5){\circle{3.0}}
 \put(83,33){$\otimes$}
 \put(85,16){\circle{3.0}}
 \put(89,20){\circle{3.0}}
 \put(100,20){\circle{3.0}}
 \put(100,35){\circle{3.0}}
 \put(72,5){\vector(1,0){11}}
 \put(72,20){\vector(1,0){15}}
 \put(70,7){\vector(0,1){11}}
 \put(85,7){\vector(0,1){7}}
 \put(91,20){\vector(1,0){7}}
 \put(85,18){\vector(0,1){15}}
 \put(87,35){\vector(1,0){11}}
 \put(100,22){\vector(0,1){11}}
 \put(65,2){$u$}
 \put(88,2){$v$}
 \put(65,21){$w$}
  \end{picture}
 \end{center}

\begin{itemize}
\item[{\bf(A8)}]
Let an $\otimes$-vertex $u$ and ordinary vertices $v,w$ be connected by
I-edges $(u,v)$ and $(v,w)$. Then $u,v,w$ belong to the small Verma
configuration (from ${\rm II}^{-1}u$ to ${\rm II}w$). Moreover, this
configuration contains exactly two $\otimes$-vertices, namely, $u$ and
${\rm I}^{-1}{\rm II}w$. (See the picture.) Symmetrically: if ordinary
vertices $u,v$ and an $\otimes$-vertex $w$ are connected by II-edges
$(u,v)$ and $(v,w)$, then $u,v,w$ belong to the small Verma
configuration, and moreover, this configuration has exactly two
$\otimes$-vertices, namely, $w$ and ${\rm II}({\rm I}^{-1}u)$.
  \end{itemize}
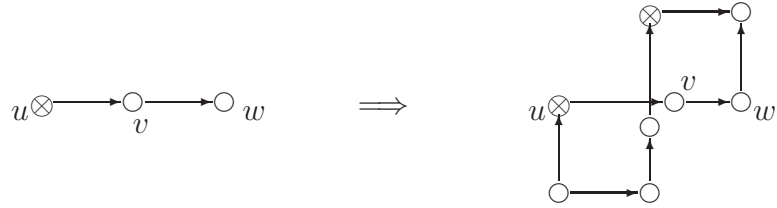
   \begin{center}
  \unitlength=0.8mm
    \begin{picture}(125,40)
 \put(3,18){$\otimes$}
 \put(35,20){\circle{3.0}}
 \put(20,20){\circle{3.0}}
 \put(7,20){\vector(1,0){11}}
 \put(22,20){\vector(1,0){11}}
 \put(0,17){$u$}
 \put(20,15){$v$}
 \put(38,17){$w$}
 \put(57,18){$\Longrightarrow$}
 \put(88,18){$\otimes$}
 \put(90,5){\circle{3.0}}
 \put(105,5){\circle{3.0}}
 \put(103,33){$\otimes$}
 \put(105,16){\circle{3.0}}
 \put(109,20){\circle{3.0}}
 \put(120,20){\circle{3.0}}
 \put(120,35){\circle{3.0}}
 \put(92,5){\vector(1,0){11}}
 \put(92,20){\vector(1,0){15}}
 \put(90,7){\vector(0,1){11}}
 \put(105,7){\vector(0,1){7}}
 \put(111,20){\vector(1,0){7}}
 \put(105,18){\vector(0,1){15}}
 \put(107,35){\vector(1,0){11}}
 \put(120,22){\vector(0,1){11}}
  \put(85,17){$u$}
 \put(110,22){$v$}
 \put(122,17){$w$}
  \end{picture}
 \end{center}

 \begin{corollary} \label{cor:A9}
Let an $\otimes$-vertex $w$ have entering I-edge $(v,w)$, and let $v$
have entering $I$-edge $(u,v)$. Then $u$ is an $\otimes$-vertex. (See
the picture.) Symmetrically: if an $\otimes$-vertex $u$ has leaving
II-edge $(u,v)$ and $v$ has leaving II-edge $(v,w)$, then $w$ is an
$\otimes$-vertex.
  \end{corollary}
   \begin{center}
  \unitlength=0.8mm
    \begin{picture}(120,10)
 \put(20,5){\circle{3.0}}
 \put(7,5){\vector(1,0){11}}
 \put(33,3){$\otimes$}
 \put(22,5){\vector(1,0){11}}
 \put(2,2){$u$}
 \put(20,0){$v$}
 \put(38,2){$w$}
 \put(55,3){$\Longrightarrow$}
 \put(77,3){$\otimes$}
 \put(95,5){\circle{3.0}}
 \put(108,3){$\otimes$}
 \put(82,5){\vector(1,0){11}}
 \put(97,5){\vector(1,0){11}}
 \put(74,1){$u$}
 \put(95,0){$v$}
 \put(113,2){$w$}
  \end{picture}
 \end{center}

Indeed (in the first claim), by (A2) for the edge $(v,w)$, $v$ has
entering II-edge $(v',v)$, and $v'$ is an $\otimes$-vertex. Then
applying the second part of~(A3) to the edge $(v',v)$, we obtain that
$u$ is an $\otimes$-vertex. The symmetric claim is shown similarly.

\begin{prop}  \label{pr:modA2}
Axioms (A0)--(A8) define precisely the set of decorated RA2-graphs.
  \end{prop}
  \begin{proof}
  To check that any decorated RA2-graph satisfies (A0)--(A8) is easy.

Conversely, let $C$ satisfy (A0)--(A8). Consider an
$\otimes$-vertex $v$ (if any). By Corollary~\ref{cor:A9}, in the
I-string containing $v$, each second vertex in backward direction
from $v$ is an $\otimes$-vertex, and in the II-string containing
$v$, each second vertex in forward direction from $v$ is an
$\otimes$-vertex (taking into account that no pair of
$\otimes$-vertices is adjacent, by (A0)(ii)). This together with
(A1)--(A5) implies that $v$ is contained in a left $A$-sail $L$
with the following properties: (a) for each vertex $u\in L$, the
maximal II-path beginning at $u$ and the maximal I-path ending at
$u$ are contained in $L$; (b) the diagonal $D(L)$ of $L$ consists
of ordinary vertices and the $\otimes$-vertices in $L$ are exactly
those that lie at odd distance from $D(L)$. Let $L$ be chosen
maximal for the given $v$. Such an $L$ does exist (with $D(L)$ not
at infinity) and is unique, which follows from the connectedness
of $C$ and (A0)(iv). Let $\Lscr$ be the set of maximal left
$A$-sails (without repetitions) constructed this way for the
$\otimes$-vertices $v$.

Next, let $C$ contain a pair of adjacent ordinary vertices.
Arguing similarly and using (A0)(iii),(A5),(A6), we extract the
set $\Rscr$ of maximal right $A$-sails whose edges are the edges
of $C$ connecting ordinary vertices.

Finally, unless $C$ consists of a single vertex (giving the trivial
RA2-graph), three cases are possible.

(i) Let $C$ have no $\otimes$-vertices. Since $C$ is connected,
$\Rscr$ consists of a unique sail $R$. Then $C=R$, and therefore,
$C$ is an RA2-graph $C(a,0)$, where $a$ is either finite or one of
$\Zset,\Zset_+,\Zset_-$ (cf. Section~\ref{sec:constr}).

(ii) Let $C$ have no pair of adjacent ordinary vertices. Then
$\Lscr$ consists of a unique sail $L$, and we have $C=R$, yielding
$C=C(0,b)$, where $b$ is either finite or one of
$\Zset,\Zset_+,\Zset_-$. Moreover, $C$ is properly decorated.

(iii) Let $C$ have both an $\otimes$-vertex and a pair of adjacent
ordinary vertices. Extend the diagonal of each sail in $\Rscr$ (in
$\Lscr$) to the corresponding path of color $\alpha$ (resp. color
$\beta$), and let $\Gamma$ be the 2-edge-colored digraph that is the
union of these paths. Applying the argument in~\cite{A-2} (in the proof
of the main structural theorem there), one shows that the small Verma
relation axioms (A7) and (A8) imply that $\Gamma$ is isomorphic to a
grid $\Gamma(a,b)$, where each of $a,b$ is either finite or one of
$\Zset,\Zset_+,\Zset_-$ (cf. Section~\ref{sec:constr}). Thus,
$C=C(a,b)$, and moreover, $C$ is properly decorated.
  \end{proof}

 \noindent
{\bf Remark 3.} One can show that Axioms~(B4),(B4$'$) follow from
Corollary~\ref{cor:A9}, and therefore, they can be excluded from the
list of axioms defining the S-graphs. In fact, property (B4) was formulated as an
axiom in order to simplify our description logically.


\section{Worm model} \label{sec:worm}

In this section we describe a model generating 2-edge-colored
digraphs; we call them {\em worm graphs}.
Vertices and edges of these graphs have a nice visualization,
which will help us to show that these graphs satisfy the axioms in
Section~\ref{sec:axiom} and that any S-graph is a worm graph. In
Appendix 1 we will take advantages from the worm model to prove
that the finite graphs among these are exactly the regular
$B_2$-crystals.

\medskip
All worm graphs are subgraphs of a universal, or {\em free}, worm graph $F$
that we now define. The vertices of $F$ are the admissible six-tuples
of integer numbers $(x',y, x''\, ;\, y',x,y'')$. Here a six-tuple is
called {\em admissible} if the following three conditions hold:

     (A) $x'$ and $x''$ are even,

     (B) $y''\ge y\ge y'$ and $x''\ge x\ge x'$,

     (C) if $y''>y$ then $x''=x$, and if $y>y'$ then $x=x'$.\medskip

It is convenient to visualize an admissible six-tuple $v=(x',y, x'';$
$y',x,y'')$ by associating to it four points in $\mathbb Z ^2$:
                        $$
        X''=(x'',y),\quad X'=(x',y),\quad Y''=(x, y'')\quad
        \text{and}\quad Y'=(x, y'),
         $$
and drawing the horizontal segment between $X'$ and $X''$ and the
vertical segment between $Y'$ and $Y''$. Then (A)--(C) are equivalent
to the following:

     1) the first coordinates of the points $X'$ and $X''$ are even;

     2) the point $X''$ lies to the right of $X'$, and the point
$Y''$ lies above $Y'$;

     3) the segments $[X',X'']$ and $[Y',Y'']$ have nonempty intersection;

     4) at least one of the following holds: $X'=X''$, $X'=Y''$, $Y'=Y''$, $Y'=X''$.

Possible cases are illustrated in the picture where $X$ stands for the
point $X'=X''$, and $Y$ for $Y'=Y''$.

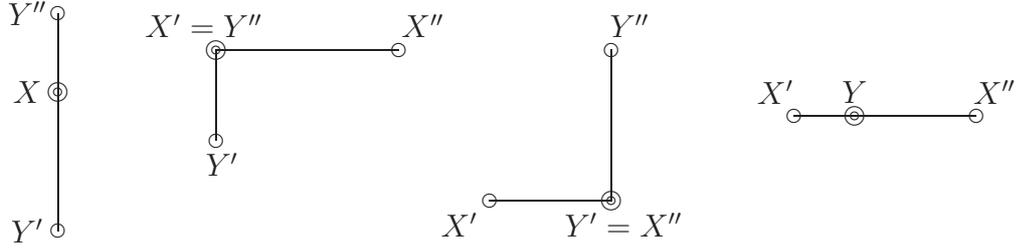
\unitlength=.800mm \special{em:linewidth 0.4pt} \linethickness{0.4pt}
 \begin{picture}(175.00,53.00)(0,5)
 \put(19.00,10.00){\circle{2.00}}
 \put(19.00,46.00){\circle{2.00}}
 \put(19.00,33.00){\circle{1.50}}
 \put(19.00,33.00){\circle{3.00}}
 \put(19,10){\line(0,1){36}}
 \put(14.00,10.00){\makebox(0,0)[cc]{$Y'$}}
 \put(14.00,46.00){\makebox(0,0)[cc]{$Y''$}}
 \put(14.00,33.00){\makebox(0,0)[cc]{$X$}}
 \put(45.00,25.00){\circle{2.00}}
 \put(45.00,40.00){\circle{1.50}}
 \put(45.00,40.00){\circle{3.00}}
 \put(75.00,40.00){\circle{2.00}}
 \put(45,25){\line(0,1){15}}
 \put(45,40){\line(1,0){30}}
 \put(46.00,21.00){\makebox(0,0)[cc]{$Y'$}}
 \put(43.00,44.00){\makebox(0,0)[cc]{$X'=Y''$}}
 \put(79.00,44.00){\makebox(0,0)[cc]{$X''$}}
 \put(90.00,15.00){\circle{2.00}}
 \put(110.00,15.00){\circle{1.50}}
 \put(110.00,15.00){\circle{3.00}}
 \put(110.00,40.00){\circle{2.00}}
 \put(90,15){\line(1,0){20}}
 \put(110,15){\line(0,1){25}}
 \put(85.00,11.00){\makebox(0,0)[cc]{$X'$}}
 \put(112.00,11.00){\makebox(0,0)[cc]{$Y'=X''$}}
 \put(113.00,44.00){\makebox(0,0)[cc]{$Y''$}}
 \put(140.00,29.00){\circle{2.00}}
 \put(150.00,29.00){\circle{1.50}}
 \put(150.00,29.00){\circle{3.00}}
 \put(170.00,29.00){\circle{2.00}}
 \put(140,29){\line(1,0){30}}
 \put(137.00,33.00){\makebox(0,0)[cc]{$X'$}}
 \put(150.00,33.00){\makebox(0,0)[cc]{$Y$}}
 \put(173.00,33.00){\makebox(0,0)[cc]{$X''$}}

 \end{picture}

 \smallskip
We call vertices of these sorts, from the left to the right in this
picture, a {\em V-worm}, a {\em VH-worm}, an {\em HV-worm}, and an {\em
H-worm}, respectively. A worm is {\em proper} if three points among
$X',X'',Y',Y''$ are different. If a worm degenerates into one point,
i.e. the corresponding six-tuple takes the form $(a,b,a;b,a,b)$, then
this vertex of $F$ is called {\em principal} (we shall see later that
such a vertex corresponds to a principal vertex in an S-graph).
\medskip

Each vertex $v=(x',y, x''\,;\,y',x,y'')$ of $F$ has two leaving edges,
colored 1 and 2, and two entering edges, colored
1 and 2 (justifying the adjective ``free''). More
precisely, the action of the operator {\bf 1} on $v$ is as follows:

\smallskip
    (i) if $2x>x'+x''$ then $x'$ increases by 2;

    (ii) if $x=x'=x''$ and $y''>y$ then $y$ increases by 1;

    (iii) otherwise $x''$ increases by 2

\smallskip
\noindent (preserving the other entries).

So in case of a proper HV-worm, the point $X'$ moves by two
positions to the right; in case of a VH-worm, the point $X''$
moves by two positions to the right; in case of a V-worm with
$X\ne Y''$, the point $X$ moves by one position up. The case of
H-worms is a bit tricky: we move (by two) that of the points
$X',X''$ which is farther from $Y$; if they are equidistant from
$Y$, then the point $X''$ moves. One can check that the operator
{\bf 1} is invertible.

\smallskip
In its turn, the action of {\bf 2} on $v$ is as follows:

\smallskip
     (iv) if $2y>y'+y''$, then $y'$ increases by 1;

     (v) if $y''=y=y'$ and $x''>x$, then $x$ increases by 1;

     (vi) otherwise $y''$ increases by 1.

\smallskip
\noindent So the operator {\bf 2} shifts $Y'$ ($Y''$) by one position
up in the proper VH-case (resp. in the HV-case) and shifts $Y$ by one
position to the right in the H-case with $Y<X''$. In the V-case, {\bf
2} shifts, by one position up, that of the points $Y',Y''$ which is
farther from $X$; if they are equidistant from $X$, then $Y''$ moves.
The operator {\bf 2} is also invertible.

\medskip
 \noindent
{\bf Remark 4.} Associate to a six-tuple $(x',y,x''\,;\,y',x,y'')$ the
pair $(x'/2+x''/2+y, y'+y''+x)$. This gives a mapping from the vertex
set of $F$ to $\Zset^2$ such that the 1-edges and 2-edges of $F$ are
congruent to the vectors $(1,0)$ and $(0,1)$, respectively; cf.
Remark~1 in Section~\ref{sec:constr}.

\medskip
Consider a string (maximal path) $P$ colored 1. One can see that
$P$ contains a V-worm or an H-worm.

Suppose $P$ contains a V-worm. When moving along this string, the
``virtual'' worm takes stages $HV$, $V$ and $VH$, in this order.
The vertical segment $[Y',Y'']$ is an invariant of the string.
Moreover, the string has a natural ``center''. When the distance
$\Vert Y'-Y''\Vert$ between $Y'$ and $Y''$ is even, this center is
defined to be the V-worm in which the double point $X$ occurs in
the middle of the vertical segment $[Y',Y'']$. When the distance
is odd, we define the center to be the edge formed by the
corresponding pair of V-worms (with $X$ lying at distance $\frac12$
below and above the middle point of $[Y',Y'']$). Thus, any
edge in $P$ is located either before the center or after the
center, unless it is the central edge itself.

If $P$ contains an H-worm, then {\em all} vertices of $P$ are H-worms
as well. An invariant of $P$ is the distance from $Y$ to the closer of
$X',X''$. Then $P$ has a natural center, the H-worm with $\Vert
X'-Y\Vert=\Vert Y-X''\Vert$ (such a worm exists since $\Vert
X'-X''\Vert$ is even).\medskip

The strings with color 2 have analogous structure and invariants.
More precisely, if a 2-string $Q$ contains an H-worm, then the segment
$[X',X'']$ is an invariant of $Q$, and the worm $(X',X'',Y=\frac
{X'+X''}2)$ is its center. If $Q$ contains a V-worm, then all vertices
of $Q$ are V-worms as well.  Then $Q$ has as an invariant the distance
from $X$ to the closer of $Y',Y''$. The center of $Q$ is defined to be
the V-worm with equal distances from $X$ to $Y'$ and to $Y''$.

Thus, the sets of central vertices for 1-strings and for 2-strings are
the same, and each central vertex is represented by a {\em symmetric}
worm. In particular, any principal vertex is central.

From the above observations we immediately obtain the following

 \begin{corollary} \label{cor:F}
The free worm graph $F$ satisfies Axioms (B0)--(B2).
  \end{corollary}

{\bf A restricted worm model: boundary conditions.}

 \medskip
We can impose natural boundary conditions on six-tuples. Let $a_1$,
$b_1$, $a_2$, $b_2$ be integers or $\pm\infty$ satisfying $a_1\le b_1$
and $a_2\le b_2$. Consider the set $W(a_1 ,b_1 ;a_2 ,b_2 )$ of
admissible six-tuples $(x',y,x''; y',x,y'')$ satisfying
  $$
2a_1 \le x', x''\le 2b_1 \quad\mbox{and}\quad a_2 \le y', y''\le b_2 .
  $$
In terms of worms, these conditions say that the worms live in the
(possibly infinite) rectangle $[2a_1 ,2b_1 ]\times [a_2 ,b_2 ]$. The
set $W(a_1 ,b_1 ;a_2 ,b_2 )$ is extended to a 2-edge-colored graph by
inducing the corresponding edges from $F$. (The operator {\bf 1} or
{\bf 2} becomes not applicable to a worm if its action would cause
trespassing the boundary of the rectangle.) One can see that the
centers of all 1- and 2-strings of $F$ intersecting the obtained graph
lie in the latter, and that Corollary~\ref{cor:F} remains valid for it.

If the boundaries are $a_1 =a_2=-\infty $ and $b_1 =b_2 =\infty $, we
have the entire graph $F$. If all $a_i$ and $b_i$ are finite, we obtain
a finite graph. If the boundaries $a_i$ and $b_i$ are shifted by the
same number, we obtain an isomorphic graph. By this reason, when both
$a_1$ and $a_2$ are finite,  it is convenient to assume that $a_1 =a_2
=0$, and we may denote the corresponding graph as $W(b_1,b_2)$.  The
graph $W=W(b_1,b _2)$ has the source (``origin'') $O=(0,0,0;0,0,0)$ (it
is easy to see that one can reach any vertex of $W$ from $O$). Also one
can see that
  \begin{itemize}
\item[($\ast$)]
the 1-string beginning at $O$ has $b_1$ edges, and the 2-string
beginning at $O$ has $b_2$ edges,
  \end{itemize}
justifying the choice of $b_1,b_2$ as the parameters of $W$.

If the numbers $b_1$ and $b_2$ are finite as well, the graph $W$ has
the sink $T=(2b_1 ,b_2 ,2b_1 ; b_2 , 2b_1 ,b_2 )$. One can see that $W$
is the interval of $F$ between the two principal vertices $O$ and $T$.

\medskip
{\em Examples of finite worm graphs.}

 \medskip
The graph $W(0,0)$ consists of a single vertex.

 \medskip
The graph $W(1,0)$ is formed by 5 worms and their transformations as
illustrated in Fig.~\ref{fig:fig6}.

\begin{figure}[htb]
 \begin{center}
\unitlength=.8mm \special{em:linewidth 0.4pt}
\linethickness{0.4pt}
    \begin{picture}(150.00,40.00)(15,5)
    {
\begin{picture}(91,30)(10,-5)
\put(30.00,5.00){\circle*{3}} \put(50.00,5.00){\circle{1.00}}
\put(50.00,5.00){\circle{2}} \put(50.00,5.00){\circle{3.00}}
\put(70.00,5.00){\circle{1.00}} \put(50.00,15.00){\circle{1.00}}
\put(70.00,15.00){\circle{1.00}} \put(60.00,15.00){\circle{1.00}}
\put(60.00,15.00){\circle{2}} \put(50.00,25.00){\circle{1.00}}
\put(70.00,25.00){\circle{1.00}} \put(70.00,25.00){\circle{2}}
\put(70.00,25.00){\circle{3.00}} \put(90.00,25.00){\circle*{3}}
\put(51,5){\line(1,0){18}}
\put(51,15){\line(1,0){18}}
\put(51,25){\line(1,0){18}}
\put(35.00,5.00){\vector(1,0){11.00}}
\put(60.00,7.00){\vector(0,1){6.00}}
\put(60.00,17.00){\vector(0,1){6.00}}
\put(75.00,25.00){\vector(1,0){11.00}}
\put(40.00,7.00){\makebox(0,0)[cc]{\bf 1}}
\put(63.00,10.00){\makebox(0,0)[cc]{\bf 2}}
\put(63.00,20.00){\makebox(0,0)[cc]{\bf 2}}
\put(80.00,22.00){\makebox(0,0)[cc]{\bf 1}}
\end{picture}
\put(10.00,15.00){$\Longrightarrow$}
\begin{picture}(106.00,43.00)(40,0)
\put(75.00,5.00){\circle*{3.00}} \put(90.00,5.00){\circle{2.00}}
\put(90.00,20.00){\circle*{2.00}} \put(90.00,35.00){\circle{2.00}}
\put(105.00,35.00){\circle*{3.00}}
\put(76.00,5.00){\vector(1,0){13.00}}
\put(90.00,6.00){\vector(0,1){13.00}}
\put(90.00,21.00){\vector(0,1){13.00}}
\put(91.00,35.00){\vector(1,0){13.00}}
\end{picture}
   }
   \end{picture}
\end{center}
\caption{$W(1,0)$}
\label{fig:fig6}
  \end{figure}

The graph $W(0,1)$ is formed by 4 worms and their transformations as
illustrated in Fig.~\ref{fig:fig7}.

\begin{figure}[htb]
 \begin{center}
\unitlength=.8mm \special{em:linewidth 0.4pt} \linethickness{0.4pt}
    \begin{picture}(150.00,45.00)(15,0)
    {
\begin{picture}(90.00,46)(20,-10)
\put(60.00,-10){\circle*{3.00}} \put(60,-6){\vector(0,1){7}}
\put(60.00,5.00){\circle{1.00}} \put(60.00,5.00){\circle{2.00}}
\put(60.00,5.00){\circle{3.00}} \put(60.00,15.00){\circle{1.00}}

\put(80.00,15.00){\circle{1.00}} \put(80.00,15.00){\circle{2.00}}
\put(80.00,15.00){\circle{3.00}} \put(80.00,5.00){\circle{1.00}}

\put(80.00,30){\circle*{3.00}} \put(80,19){\vector(0,1){7}}

\put(63.00,10.00){\vector(1,0){14.00}}

\put(63,-3){\makebox(0,0)[cc]{\bf 2}}
\put(70.00,13.00){\makebox(0,0)[cc]{\bf 1}}
\put(83.00,22.00){\makebox(0,0)[cc]{\bf 2}}
\bezier{20}(60.00,6.00)(60.00,10.00)(60.00,14.00)
\bezier{20}(80.00,6.00)(80.00,10.00)(80.00,14.00)
\end{picture}
\put(0.00,15.00){$\Longrightarrow$}
\begin{picture}(60.00,42.00)(10,0)
\put(40.00,5.00){\circle*{3.00}} \put(40.00,20.00){\circle{2.00}}
\put(55.00,20.00){\circle{2.00}} \put(55.00,35.00){\circle*{3.00}}
\put(40.00,6.00){\vector(0,1){13.00}}
\put(48.00,20.00){\makebox(0,0)[cc]{\scriptsize $\blacklozenge$}}
\put(41.00,20.00){\vector(1,0){13.00}}
\put(55.00,21.00){\vector(0,1){13.00}}
\end{picture}
   }
   \end{picture}
\end{center}
\caption{$W(0,1)$}
\label{fig:fig7}
  \end{figure}

Therefore, $W(1,0)$ and $W(0,1)$ are the same as the sail graphs
$S(1,0)$ and $S(0,1)$, respectively (cf. Fig.~\ref{fig:fig4}). A more tiresome, but
useful, exercise is to construct the worm graph $W(1,1)$ and check that
it is equal to the graph $S(1,1)$.

 \medskip
 \noindent
{\bf Remark 5.} For an edge $e=(u,v)$ of a finite worm-graph, it is not
difficult to compute the numbers $\Delta t(e)$ and $\Delta h(e)$
(defined in Remark~2), which depend on the color of $e$ and the form of
the worm $u$. More precisely, when $e$ has color 1: (i) $\Delta
t(e)=0$ and $\Delta h(e)=2$ if $u$ is a VH-worm or a V-worm with $\Vert
X-Y'\Vert\ge \Vert Y''-X\Vert$ or an H-worm with $\Vert Y-X'\Vert \le
\Vert X''-Y\Vert$; (ii) $\Delta t(e)=-1$ and $\Delta h(e)=1$ if $u$ is
a V-worm with $\Vert X-Y'\Vert =\Vert Y''-X\Vert-1$; and (iii) $\Delta
t(e)=-2$ and $\Delta h(e)=0$ in the other cases. When $e$ has color
2: (iv) $\Delta t(e)=-1$ and $\Delta h(e)=0$ if $u$ is a VH-worm
or a V-worm with $\Vert X-Y'\Vert\ge \Vert Y''-X\Vert$ or an H-worm
with $\Vert Y-X'\Vert \le \Vert X''-Y\Vert$; and (v) $\Delta t(e)=0$
and $\Delta h(e)=1$ in the other cases.

 \medskip
Next we show the following key property.

\begin{prop} \label{pr:W-ax}
Any worm graph $W$ satisfies Axioms (B3),(B4),(B3$'$),(B4$'$),(BA).
  \end{prop}
  \begin{proof}
We denote the quadruple corresponding to a vertex (worm) $v$ of $W$ by
$q(v)=(X'(v),X''(v),Y'(v),Y''(v))$.

First we verify Axiom~(B3). Consider a left 1-edge $(u,v)$ in $W$.
Since the vertex $u$ is left and the edge $(u,v)$ is not central, only
three cases are possible:

(a) $u$ is a V-worm, and $\Vert Y''(u)-X(u)\Vert\ge \Vert
X(u)-Y'(u)\Vert+2$;

(b) $u$ is a proper HV-worm;

(c) $u$ is an H-worm, and $\Vert X'(u)-Y(u)\Vert\ge \Vert
Y(u)-X''(u)\Vert+2$.

\smallskip
 \noindent
The quadruple $q(v)$ is obtained from $q(u)$ by shifting $X$ by 1 up in
case (a), and shifting $X'$ by 2 to the right in cases (b),(c)
(preserving the other entries). Suppose $u$ has leaving 2-edge
$(u,u')$. Then $q(u')$ is obtained from $q(u)$ by shifting $Y''$ by 1
up in cases (a),(b) and in the subcase of (c) with $Y(u)=X''(u)$, and
by shifting $Y$ by 1 to the right in the subcase of (c) with $Y(u)\ne
X''(u)$. Form the quadruple $q=(X'(v),X''(v),Y'(u'),Y''(u'))$. It is
straightforward to check that in all cases $q$ determines a feasible
worm $w$, and that $w={\bf 2}v={\bf 1}u'$, as required.

Now suppose that $v$ has leaving 2-edge $(v,v')$. Then $q(v')$ is
obtained from $q(v)$ by shifting $Y''$ by 1 up in cases (a),(b) and in
the subcase of (c) with $Y(v)=X''(v)$, and by shifting $Y$ by 1 to the
right in the subcase of (c) with $Y(v)\ne X''(v)$. Again, one easily
checks that in all cases the quadruple $(X'(u),X''(u),Y'(v'),Y''(v'))$
determines a feasible worm $u'$, that $u'={\bf 2}u$, and that $v'={\bf
1}u'$, as required.

This gives part (i) of Axiom~(B3). Part (ii) of this axiom is shown in
a similar way.

Validity of the dual axiom (B3$'$) follows from (B3) and the symmetry of
$F$ (in the sense that reversing the edges makes the graph isomorphic
to $F$). To verify (B4),(B4$'$) is not necessary; see Remark~3.

\smallskip
Next we verify Axiom~(BA). As is said above, the central vertices $v$ of $W$
are defined to be the symmetric V-worms (i.e., with $\Vert Y'-Y''\Vert$
even and with $X$ in the middle of $[Y',Y'']$) and the symmetric
H-worms (i.e., with $Y$ in the middle of $[X',X'']$). Such a $v$ is an ordinary
vertex of the central graph $C(W)$ of $W$, and we identify it with
the corresponding even-length interval $J(v)=[U(v),V(v)]$ of the form $[Y',Y'']$ or
$[X',X'']$. The central 1-edges of $W$ are the pairs $(u,v)$ of
V-worms with the same $Y'$ and the same $Y''$ and such that $\Vert
Y'-Y''\Vert$ is odd and $\Vert Y'-X(u)\Vert=\Vert Y'-X(v)\Vert-1=\Vert
X(v)-Y''\Vert$. Such an edge $e$ generates an $\otimes$-vertex of $C(W)$,
which we now denote by $e$ as well, and we identify it with the
odd-length vertical interval $[Y',Y'']$ denoted by $J(e)=[U(e),V(e)]$.

Thus, we have an {\em interval model} (simplifying the worm model) to
represent the vertices of $C(W)$. The $\otimes$-vertices of
$C(W)$ are the odd-length vertical intervals, or, briefly, the {\em
odd} intervals, and the ordinary vertices are the even-length (vertical
or horizontal) intervals, or the {\em even} intervals, in the model.
The edges of $C(W)$ correspond to the following transformations of
intervals:
 \begin{itemize}
 \item[(E1)]
for a non-degenerate horizontal interval $J=[U,V]$, the operator I
(when applicable) increases $J$ by shifting $V$ by 2 to the right,
while II decreases $J$ by shifting $U$ by 2 to the right;
  \item[(E2)]
for a non-degenerate vertical interval $J=[U,V]$, the operator I
decreases $J$ by shifting $U$ by 1 up, while II (when applicable)
increases $J$ by shifting $V$ by 1 up;
  \item[(E3)]
 for a degenerate interval $J=[U,V]$ (i.e., $U=V$), the operator I
(when applicable) shifts $V$ by 2 to the right, while II (when
applicable) shifts $V$ by 1 up.
  \end{itemize}

This can be seen from the following observations. For a non-degenerate
central H-worm $v$, the interval ${\rm I}(J(v))$ is obtained when we
apply to $v$ the operator {\bf 1} followed by {\bf 2} (and similarly if
$v$ is degenerate), and the interval ${\rm II}(J(v))$ is obtained when
we apply to $v$ the operator {\bf 2} followed by {\bf 1}. For a
non-degenerate central V-worm $v$, the interval ${\rm I}(J(v))$
corresponds to the central 1-edge with the beginning ${\bf 1}^{-1}{\bf
2}{\bf 1}v$ and the end ${\bf 2}{\bf 1}v$, and the interval ${\rm
II}(J(v))$ corresponds to the central 1-edge with the beginning ${\bf
2}v$ and the end ${\bf 1}{\bf 2}v$ (and similarly if $v$ is
degenerate). And for a central 1-edge $e=(u,v)$, the interval ${\rm
I}(J(e))$ corresponds to the central V-worm ${\bf 2}v$, while the
interval ${\rm II}(J(e))$ corresponds to the central V-worm ${\bf
1}{\bf 2}u$.

We have to show that $C(W)$ satisfy (A0)--(A8). Property (A0) is easy.
As is seen from (E1)--(E3), the central graph $C(F)$ of the free worm
graph $F$ remains essentially the same when we reverse the edges and
replace the operator I by ${\rm II}^{-1}$, and II by ${\rm I}^{-1}$.
Therefore, in (A2)--(A8), it suffices to verify only the first parts of
these axioms.

 \smallskip
(i) If $J=[U,V]$ is an odd vertical interval, then both I and ${\rm
II}^{-1}$ are applicable to $J$ (the former lifts $U$ by 1 and the
latter descends $V$ by 1). Also the resulting vertical intervals are even. This
gives~(A1).

 \smallskip
(ii) Let $J=[U,V]$ be an odd vertical interval and ${\rm I}^{-1}$ is
applicable to $J$ (i.e., $U$ is not on the bottom of the rectangle
bounding $W$). Let $U'=U-(0,1)$ and $V'=V-(0,1)$. Then ${\rm
I}^{-1}J=[U',V]$, ${\rm II}^{-1}J=[U,V']$ and ${\rm II}^{-1}[U',V]
={\rm I}^{-1}[U,V']=[U',V']$. Also $[U',V']$ is odd. This gives the
first part of (A2).

 \smallskip
(iii) Let $J=[U,V]$ and $J'=[U',V]$ be odd and even vertical intervals
with $U'=U+(0,1)$. If II is applicable to $J$ or $J'$, then the point
$V'=V+(0,1)$ belongs to the rectangle, and the first part of~(A3)
follows by arguing as in (ii). The first part of~(A4) is shown in a
similar way.

\smallskip
(iv) Let $J=[U,V]$ and $J'=[U',V']$ be even intervals connected by
I-edge $(J,J')$. Then $J,J'$ are horizontal intervals, $U=U'$, and
$V'=V+(2,0)$. The feasible and even interval $[U+(2,0),V']$ is just
$II(J')$, yielding the first part of~(A5).

 \smallskip
(v) The first part of (A6) is shown by a method as in (ii), with the
difference that now we deal with horizontal intervals.

 \smallskip
(vi) Let $J,J',J''$ be the intervals for $u,v,w$ (respectively) as in
the first part of~(A7). Then $J,J'$ are even and $J''$ is odd. Moreover,
since these intervals are connected by I-edge $(J,J')$ and by II-edge
$(J,J'')$, the only possible case is that $J$ is a degenerate interval
$[U,U]$, and therefore, $J'$ and $J''$ are the intervals $[U,V]$ and
$[U,V']$, respectively, where $V=U+(2,0)$ and $V'=U+(0,1)$. In the
rectangle bounding $W$ take the point $Z$ such that $[V,Z]$ is a
vertical interval and $[V',Z]$ is a horizontal interval; see the
picture. We have: $[V,V]=II(J')$, $[V,Z]=II[V,V]$, $[V',V']=I(J'')$,
$[V',Z]=I[V',V']$ and $[Z,Z]=I[V,Z]=II[V',Z]$. Also among the eight
intervals above, only $J''$ and $[V,Z]$ are odd. This implies the first
part of~(A7).

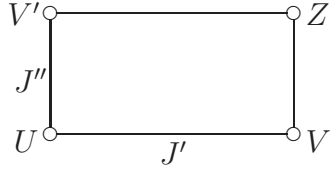
\unitlength=.8mm \special{em:linewidth 0.4pt} \linethickness{0.4pt}
\begin{picture}(60,32)(-50,0)
  \put(10,5){\circle{2.00}}
  \put(50,5){\circle{2.00}}
 \put(10,25){\circle{2.00}}
 \put(50,25){\circle{2.00}}
 \put(11,5){\line(1,0){38.00}}
 \put(11,25){\line(1,0){38.00}}
 \put(10,6){\line(0,1){18.00}}
 \put(50,6){\line(0,1){18.00}}
 \put(4,2.00){$U$}
 \put(4,12){$J''$}
 \put(3,23){$V'$}
 \put(52,2){$V$}
 \put(28,0){$J'$}
 \put(52,23){$Z$}
 \end{picture}

 \smallskip
(vii) Let $J,J',J''$ be the intervals for $u,v,w$ (respectively) as in
the first part of~(A8). Then $J$ is odd and $J',J''$ are even. Since
these intervals are connected by I-edges $(J,J')$ and $(J',J'')$, $J'$
is a degenerate interval $[V,V]$, and therefore, $J$ and $J''$ are the
intervals $[U,V]$ and $[V,V']$, respectively, where $U=V-(0,1)$ and
$V'=V+(2,0)$. Take the point $Z$ such that $[U,Z]$ is a horizontal
interval and $[Z,V']$ is a vertical interval. Then the eight intervals
(including degenerate ones) in the rectangle $UVV'Z$ give the small
Verma configuration as required in the first part of (A8).

 \smallskip
This completes the proof of the proposition.
  \end{proof}

\begin{theorem}  \label{tm:S-W}
The sets of S-graphs and worm graphs are the same.
  \end{theorem}
  \begin{proof}
Theorem~\ref{tm:ax-S}, Proposition~\ref{pr:W-ax} and
Corollary~\ref{cor:F} imply that any worm graph is an S-graph. To see
that the sets of worm graphs and S-graphs coincide, we first observe this
coincidence for the case of finite graphs. Indeed, each finite S-graph
$S$ is uniquely determined by the pair of its parameters
$a,b\in\Zset_+$ (i.e., $S$ is $S(a,b)$ defined in
Section~\ref{sec:constr}). We have seen that the numbers $a,b$ are just
the lengths of 1-string and 2-string, respectively, from the source of
$S$. In view of property~($\ast$) above, a similar behavior takes place for the finite
worm graphs $W(a,b)$ (where $(a,b)$ runs over $\Zset_+\times\Zset_+$). So the sets of finite S-graphs and worm graphs coincide.

For an infinite S-graph $S$, one can arrange an infinite sequence
$S_1,S_2,\ldots$ of intervals between principal vertices of $S$ (see
Remark~1) such that: (i) each $S_i$ is isomorphic to a finite S-graph
$S(a_i,b_i)$ and contains $S_{i-1}$ as an interval between its
principal vertices, and (ii) this sequence tends to $S$. We can arrange
a corresponding sequence $W_1\subset W_2\subset \ldots$ of restricted
finite worm subgraphs of the free graph $F$, where each $W_i$ is a
shift of $W(a_i,b_i)$. Then this sequence tends to a worm graph $W$
which is isomorphic to $S$, and the result follows.
  \end{proof}

 \medskip
 \noindent
{\bf Remark 6.}  In fact, the restricted worm model gives an alternative local axiomatics for finite $B_2$-crystals. Indeed, in the worm model each vertex $v$ of a graph is endowed with a six-tuple $\tau(v)$ of integers $(x,x',x'',y,y',y'')$ satisfying conditions ~(A),(B),(C). Local conditions in (i)--(vi) prescribe how the six-tuple $\tau(v)$ can change under the action of crystal operators (thus defining the edges of colors 1 and 2 entering and leaving $v$). 

\section*{Appendix 1: A relation to Littelmann's cones}

We use the brief notation $W(\infty)$ for the restricted worm graph
$W(0,\infty\,;\,0,\infty)$ (see Section~\ref{sec:worm}). This graph has
the following properties: it admits a weight mapping (see Remark~4); it
has one minimal vertex (the source), which corresponds to the ``origin''
six-tuple $O=(0,0,0\,;\,0,0,0)$; each of its vertices has two leaving
edges, of color 1 and color 2. Also from the worm model it
follows that the finite worm graphs (considered up to isomorphism) are
parameterized by the pairs $(a,b)\in\Zset_+^2$, and each $W=W(a,b)$
satisfies: (a) the lengths of the 1-string and 2-string beginning at $O$
are equal to $a$ and $b$, respectively; (b) $W$ is the interval of
$W(\infty)$ between $O$ and the principal vertex (the sink) of the form
$({\bf 212})^b({\bf 1221})^aO$ (cf. the fundamental graphs $W(1,0)$ and
$W(0,1)$); (c) for each edge $e$ of $W$, the number $\Delta t(e)-\Delta
h(e)$ (see Remarks~2 and~5) is as prescribed by the Cartan matrix $B_2$.

In light of the above properties, in order to establish that the finite
worm graphs (vis. S-graphs) are precisely the regular $B_2$-crystals, it
suffices to show the following.

  \begin{theorem} \label{tm:W-B2}
$W(\infty)$ is isomorphic to the graph $B(\infty)$ for $B_2$ type.
  \end{theorem}

(The graph $B(\infty)$, the inductive limit of $B_2$-crystal bases
under their inclusions agreeable with growing the dominant
weights, has similar properties as those for $W(\infty)$ above:
$B(\infty)$ has one source, admits a weight mapping and has
leaving edges of both colors at each vertex, and each regular
$B_2$-crystal is an interval of $B(\infty)$ beginning at the
source (cf.~\cite[Ch.7]{Kbook}).)

Our aim is to prove that $W(\infty)$ satisfies Littelmann's
characterization of $B(\infty)$ for $B_2$ type, formulated in terms of two
cones in $\Zset^4$.

Let $G=(V_G,E_G)$ be a 2-edge colored directed graph, with colors 1
and 2, which admits a weight mapping and such that each
monochromatic component of $G$ is a path having a beginning vertex; for
convenience we call such a graph {\em normal}. Let $\tone$ ($\ttwo$)
denote the operator which brings each vertex $v$ of $G$ to the beginning
of the 1-string (resp. 2-string) containing $v$. Define the numbers
 \begin{gather*}
 a_1(v):=t_2(v), \quad a_2(v):=t_1(\ttwo v), \quad a_3(v):=t_2(\tone\ttwo v),
 \quad  a_4(v):=t_1(\ttwo\tone\ttwo v);  \\
 b_1(v):=t_1(v), \quad b_2(v):=t_2(\tone v), \quad b_3(v):=t_1(\ttwo\tone v),
 \quad b_4(v):=t_2(\tone\ttwo\tone v),
  \end{gather*}
where for a vertex $u$, $t_i(u)$ is the length of the $i$-colored path
from $\tilde i(u)$ to $u$. We denote the quadruple
$(a_1(v),a_2(v),a_3(v),a_4(v))$ by ${\bf a}(v)$, and the quadruple
$(b_1(v),b_2(v),b_3(v),b_4(v))$ by ${\bf b}(v)$.

  \begin{theorem}[\cite{Littl}] \label{tm:litt}
Let $G=(V_G,E_G)$ be a normal graph with one source $O$ which possesses
the following properties:

{\rm (L1)} For each vertex $v$, both $\tone\ttwo\tone\ttwo v$ and
$\ttwo\tone\ttwo\tone v$ coincide with the source $O$; in particular,
$v$ is represented as
  $$
v={\bf 2}^{a _1(v)}{\bf 1}^{a_2(v)}{\bf 2}^{a_3(v)}{\bf 1}^{a_4(v)}O=
 {\bf 1}^{b _1(v)}{\bf 2}^{b_2(v)}{\bf 1}^{b_3(v)}{\bf 2}^{b_4(v)}O.
         $$

{\rm (L2)} The set $C_1:=\{{\bf a}(v)\colon v\in V_G\}$ consists of the
quadruples $(a_1,a_2,a_3,a_4)$ of nonnegative integers satisfying
   \begin{equation} \label{eq:C1}
   2a_2\ge a_3\ge 2a_4.
   \end{equation}

{\rm(L3)} The set $C_2:=\{{\bf b}(v)\colon v\in V_G\}$ consists of the
quadruples $(b_1,b_2,b_3,b_4)$ of nonnegative integers satisfying
   \begin{equation} \label{eq:C2}
   b_2\ge b_3\ge b_4.
   \end{equation}

{\rm(L4)} For each vertex $v$, the vectors ${\bf
a}(v)=(a_1,a_2,a_3,a_4)$ and ${\bf b}(v)=(b_1,b_2,b_3,b_4)$ satisfy the
following relations:
   \begin{equation} \label{eq:rel}
          b_3 =\min(a_2 ,2a_2 -a_3 +a_4, a_1+a_4 ) \quad\mbox{and}
  \quad b_4=\min(a_1 ,2a_2 -a _3,a_3 -2a_4).
   \end{equation}

Then $G$ is isomorphic to the graph $B(\infty)$ for $B_2$ type.
Conversely, $B(\infty)$ satisfies (L1)--(L4).
  \end{theorem}

Note that (L1) implies that all quadruples ${\bf a}(v)$, $v\in V_G$,
are different, and similarly for the quadruples ${\bf b}(v)$. Also
since $G$ has a weight mapping, for each vertex $v$, we have $a_1(v)
+a_3(v) =b_2(v) +b_4(v)$ and $a_2(v)+a_4(v)=b_1(v)+b_3(v)$, and
therefore, one can transform~(\ref{eq:rel}) into relations for $b_1$
and $b_2$.

Thus, we have to show that $W(\infty)$ satisfies conditions (L1)--(L4) in
Theorem~\ref{tm:litt}. In fact, we can reduce verification of (L2) (resp.
(L3)) to merely checking that the quadruples $\bfa(v)$ (resp. $\bfb(v)$)
of all vertices $v$ of $W(\infty)$ belong to the cone defined
by~(\ref{eq:C1}) (resp. by~(\ref{eq:C2})). Indeed, by
Theorem~\ref{tm:litt}, the vertices of $B(\infty)$ determine a bijection
$\gamma: C_1\to C_2$, and this bijection is given by~(\ref{eq:rel}).
Suppose (under validity of~(\ref{eq:rel}) for $W(\infty)$) that some
$\bfa\in C_1$ is not realized by a vertex of $W(\infty)$ (equivalently:
$\gamma(\bfa)\in C_2$ is not realized by a vertex of $W(\infty)$). Among
such elements, choose $\bfa$ for which the corresponding vertex $v$ of
$B(\infty)$ has minimum distance from the source. Clearly this distance is
nonzero, so $v$ has an entering edge $(u,v)$. Then $\bfa(u)=\bfa(u')$ for
some vertex $u'$ of $W(\infty)$. Letting for definiteness that $(u,v)$ has
color 2, take the edge $(u',v')$ with color 2 in $W(\infty)$.
Since both $\bfa(v)$ and $\bfa(v')$ are obtained from $\bfa(u)=\bfa(u')$
by increasing the first entry by one, we have $\bfa(v')=\bfa(v)=\bfa$; a
contradiction.

  \smallskip
We prove (L1) and \refeq{C1}--\refeq{rel}  for the vertices of
$W(\infty)$ by considering six possible cases of a worm $v$. We will
use the following conventions. In the transformations of $v$ below,
$X',X'',Y',Y''$ denote the points for the {\em current} worm, while the
numbers $x',y,x'',y',x,y''$ concern the {\em original} $v$. Everywhere
$a_i$ and $b_i$ stand for $a_i(v)$ and $b_i(v)$, respectively.  We
define (cf.~\refeq{rel}):
   $$
   p:=2a_2 -a_3 +a_4, \quad q:=a_1+a_4, \quad r:=2a_2 -a _3, \quad
   s:=a_3 -2a_4,
   $$
and define:
 \begin{gather*}
 v_1:=\ttwo v, \quad v_2:=\tone v_1, \quad v_3:=\ttwo v_2, \quad
 v_4:=\tone v_3; \\
  v'_1:=\tone v, \quad v'_2:=\ttwo v'_1, \quad v'_3:=\tone v'_2,
 \quad  v'_4:=\ttwo v'_3.
  \end{gather*}

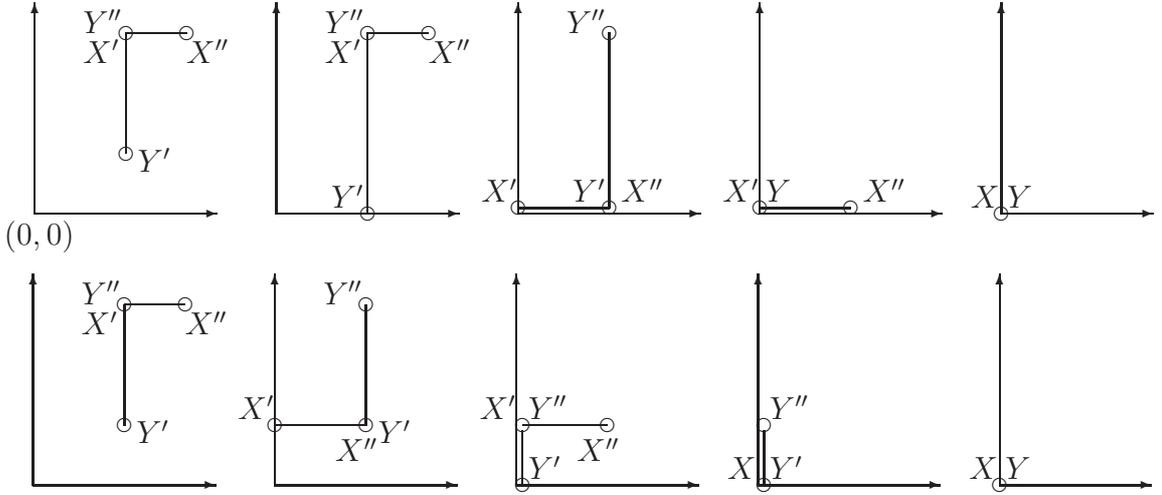
\begin{figure}[htb]
 \begin{center}
\unitlength=.8mm \special{em:linewidth 0.4pt} \linethickness{0.4pt}
 \begin{picture}(180,90)
  {
  \begin{picture}(180,45)(5,-45)
  {
 \begin{picture}(30,40)(5,0)    
 \put(5,5){\vector(1,0){30}}
 \put(5,5){\vector(0,1){35}}
 \put(20,15){\circle{2}}
 \put(20,35){\circle{2}}
 \put(30,35){\circle{2}}
 \put(20,15){\line(0,1){20}}
 \put(20,35){\line(1,0){10}}
 \put(0,0){$(0,0)$}
 \put(22,12){$Y'$}
 \put(13,35){$Y''$}
 \put(13,30){$X'$}
 \put(30,30){$X''$}
  \end{picture}
 \begin{picture}(30,40)(-3,0)    
 \put(5,5){\vector(1,0){30}}
 \put(5,5){\vector(0,1){35}}
 \put(20,5){\circle{2}}
 \put(20,35){\circle{2}}
 \put(30,35){\circle{2}}
 \put(20,5){\line(0,1){30}}
 \put(20,35){\line(1,0){10}}
 \put(14,6){$Y'$}
 \put(13,35){$Y''$}
 \put(13,30){$X'$}
 \put(30,30){$X''$}
    \end{picture}
 \begin{picture}(30,40)(-11,0)    
 \put(5,5){\vector(1,0){30}}
 \put(5,5){\vector(0,1){35}}
 \put(20,6){\circle{2}}
 \put(20,35){\circle{2}}
 \put(20,6){\line(0,1){29}}
 \put(5,6){\circle{2}}
 \put(5,6){\line(1,0){15}}
 \put(22,7){$X''$}
 \put(13,35){$Y''$}
 \put(14,7){$Y'$}
 \put(-1,7){$X'$}
   \end{picture}
 \begin{picture}(30,40)(-19,0)     
 \put(5,5){\vector(1,0){30}}
 \put(5,5){\vector(0,1){35}}
 \put(20,6){\circle{2}}
 \put(5,6){\circle{2}}
 \put(5,6){\line(1,0){15}}
 \put(22,7){$X''$}
 \put(6,7){$Y$}
 \put(-1,7){$X'$}
   \end{picture}
 \begin{picture}(30,40)(-27,0)      
 \put(5,5){\vector(1,0){25}}
 \put(5,5){\vector(0,1){35}}
 \put(5,5){\circle{2}}
 \put(6,6){$Y$}
 \put(0,6){$X$}
   \end{picture}
   }
   \end{picture}
  \begin{picture}(180,45)(187,0)
  {
 \begin{picture}(30,40)(5,0)       
 \put(5,5){\vector(1,0){30}}
 \put(5,5){\vector(0,1){35}}
 \put(20,15){\circle{2}}
 \put(20,35){\circle{2}}
 \put(30,35){\circle{2}}
 \put(20,15){\line(0,1){20}}
 \put(20,35){\line(1,0){10}}
 \put(22,12){$Y'$}
 \put(13,35){$Y''$}
 \put(13,30){$X'$}
 \put(30,30){$X''$}
  \end{picture}
 \begin{picture}(30,40)(-3,0)        
 \put(5,5){\vector(1,0){30}}
 \put(5,5){\vector(0,1){35}}
 \put(5,15){\circle{2}}
 \put(20,35){\circle{2}}
 \put(20,15){\circle{2}}
 \put(20,15){\line(0,1){20}}
 \put(5,15){\line(1,0){15}}
 \put(22,12){$Y'$}
 \put(13,35){$Y''$}
 \put(-1,16){$X'$}
 \put(15,9){$X''$}
   \end{picture}
 \begin{picture}(30,40)(-11,0)      
 \put(5,5){\vector(1,0){30}}
 \put(5,5){\vector(0,1){35}}
 \put(20,15){\circle{2}}
 \put(6,15){\circle{2}}
 \put(6,5){\circle{2}}
 \put(6,15){\line(1,0){14}}
 \put(6,5){\line(0,1){9}}
 \put(7,6){$Y'$}
 \put(7,16){$Y''$}
 \put(-1,16){$X'$}
 \put(15,9){$X''$}
   \end{picture}
 \begin{picture}(30,40)(-19,0)       
 \put(5,5){\vector(1,0){30}}
 \put(5,5){\vector(0,1){35}}
 \put(6,15){\circle{2}}
 \put(6,5){\circle{2}}
 \put(6,5){\line(0,1){9}}
 \put(7,6){$Y'$}
 \put(7,16){$Y''$}
 \put(0,6){$X$}
      \end{picture}
 \begin{picture}(30,40)(-27,0)       
 \put(5,5){\vector(1,0){25}}
 \put(5,5){\vector(0,1){35}}
 \put(5,5){\circle{2}}
 \put(6,6){$Y$}
 \put(0,6){$X$}
   \end{picture}
   }
  \end{picture}
 }
   \end{picture}
 \end{center}
 \caption{The transformations for a VH-worm $v$. In the upper line:
 $v\to v_1\to v_2\to v_3\to v_4=O$. In the lower line: $v\to
 v'_1\to v'_2\to v'_3\to v'_4=O$.}
\label{fig:fig8}
  \end{figure}

{\bf Case 1}: $v$ is a proper VH-worm. The chain of transformations
$v\to v_1\to v_2\to v_3\to v_4$ is illustrated in the upper line of
Fig.~\ref{fig:fig8}. More precisely, the operator $\ttwo$ applied to $v$ moves
$Y'=(x',y')$ to the point $(x',0)$. The action of $\tone$ at $v_1$
moves $X''=(x'',y'')$ to the point $(x',y'')$ (using $\frac{x''-x'}{2}$
steps), then moves the double point $X=(x',y'')$ to $(x',0)$ ($=Y'$),
and eventually moves $X'=(x',0)$ to the origin $(0,0)$ (using
$\frac{x'}{2}$ steps). The action of $\ttwo$ at $v_2$ moves
$Y''=(x',y'')$ to $(x',0)$ and then moves the double point $Y$ to
$(0,0)$. Finally, the action of $\tone$ at $v_3$ moves $X''=(x',0)$ to
$(0,0)$ (using $\frac{x'}{2}$ steps). So $v_4$ is the source $O$, as
required in~(L1). A direct count gives
  $$
  a_1=y', \quad a_2=\frac{x''}{2}+y'',\quad a_3=x'+y'', \quad
  a_4=\frac{x'}{2}.
  $$
Relation~\refeq{C1} turns into $x''+2y''\ge x'+y''\ge x'$, which is
valid because $x''\ge x'$.

In the second chain of transformations (see the lower line of
Fig.~\ref{fig:fig8}), the action of $\tone$ at $v$ moves $X''$ to
$(x',y'')$ (using $\frac{x''-x'}{2}$ steps), then moves the double point
$X$ to $(x',y')$ ($=Y'$), and then moves $X'$ to $(0,y')$ (using
$\frac{x'}{2}$ steps). The action of $\ttwo$ at $v'_1$ moves $Y''$ to
$(x',y')$, then moves the double point $Y$ to $(0,y')$, and eventually
moves $Y'$ to $(0,0)$. The action of $\tone$ at $v'_2$ moves $X''$ to
$(0,y')$ (using $\frac{x'}{2}$ steps), and then moves $X$ to $(0,0)$. And
the action of $\ttwo$ at $v'_3$ moves $Y''$ to $(0,0)$. So $v'_4=O$, as
required. We have
  $$
  b_1=\frac{x''}{2}+y''-y',\quad b_2=x'+y'',\quad b_3=\frac{x'}{2}+y',
  \quad b_4=y'.
  $$
Then~\refeq{C2} turns into $x'+y''\ge \frac{x'}{2}+y'\ge y'$, which is
valid.

Finally, we obtain
$p=(x''+2y'')-(x'+y'')+\frac{x'}{2}=x''-\frac{x'}{2}+y''$ and
$q=y'+\frac{x'}{2}$. Then $b_3=\frac{x'}{2}+y'=q$ and $b_3\le a_2,p$,
yielding the first equality in~\refeq{rel}. We also have
$r=(x''+2y'')-(x'+y'')=x''-x'+y''\ge y'$ and $s=(x'+y'')-x'=y''$. Then
$b_4=y'=a_1\le r,s$, yielding the second equality in~\refeq{rel}.

 \medskip
{\bf Case 2}: $v$ is a proper HV-worm. The transformations in this case
are examined straightforwardly as well. The action of $\ttwo$ at $v$
moves $Y''$ to $(x'',y')$, followed by moving $Y$ to $(x',y')$,
followed by moving $Y'$ to $(x',0)$. The action of $\tone$ at $v_1$
moves $X''$ to $(x',y')$, followed by moving $X$ to $(x',0)$, followed
by moving $X'$ to $(0,0)$. The action of $\ttwo$ at $v_2$ moves $Y''$
to $(x',0)$, followed by moving $Y$ to $(0,0)$. And the action of
$\tone$ at $v_3$ moves $X''$ to $(0,0)$.

In the second chain, the action of $\tone$ at $v$ moves $X'$ to
$(0,y')$. The action of $\ttwo$ at $v'_1$ moves $Y''$ to $(x'',y')$,
followed by moving $Y$ to $(0,y')$, followed by moving $Y'$ to $(0,0)$.
The action of $\tone$ at $v'_2$ moves $X''$ to $(0,y')$, followed by
moving $X$ to $(0,0)$. And the action of $\ttwo$ at $v'_3$ moves $Y''$
to $(0,0)$.

This gives $v_4=v'_4=O$ and:
   \begin{gather*}
  a_1=x''-x'+y'', \quad a_2=\frac{x''}{2}+y',\quad a_3=x'+y', \quad
  a_4=\frac{x'}{2}; \\
  b_1=\frac{x'}{2},\quad b_2=x''+y'',\quad b_3=\frac{x''}{2}+y',
  \quad b_4=y'.
  \end {gather*}

Relation~\refeq{C1} turns into $x''+2y'\ge x'+y'\ge x'$, and~\refeq{C2}
turns into $x''+y''\ge \frac{x''}{2}+y'\ge y'$, which are valid.

Finally, we have $a_2=b_3$,
$p=(x''+2y')-(x'+y')+\frac{x'}{2}=x''-\frac{x'}{2}+y'\ge b_3$, and
$q=(x''-x'+y'')+\frac{x'}{2}\ge b_3$, yielding the first equality
in~\refeq{rel}. Also $a_1=x''-x'+y''\ge y'=b_4$, $r=(x''+2y')-(x'+y')\ge
x''-x'+y'\ge b_4$, and $s=(x'+y')-x'=y'=b_4$, yielding the second equality
in~\refeq{rel}.

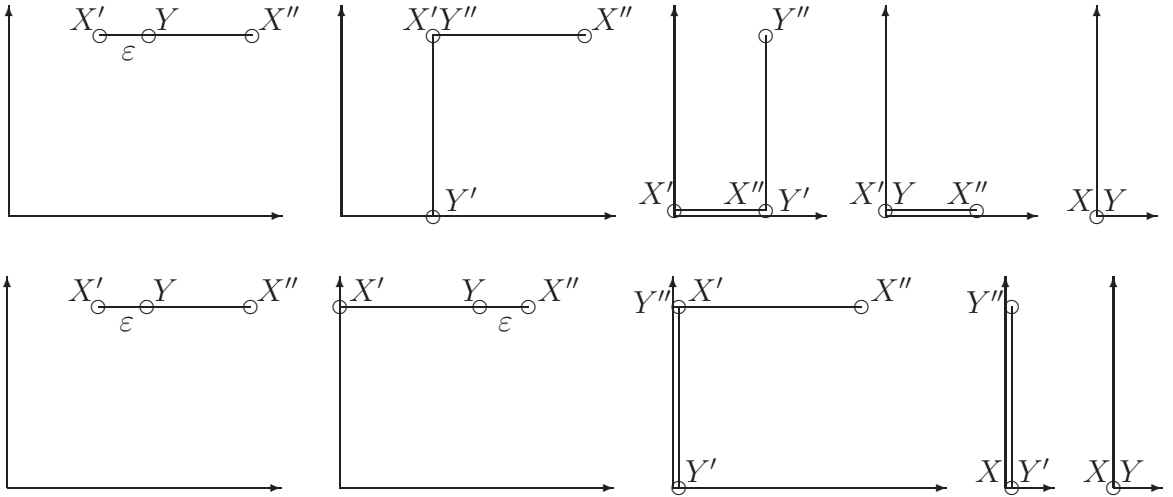
\begin{figure}[htb]
 \begin{center}
\unitlength=.8mm \special{em:linewidth 0.4pt} \linethickness{0.4pt}
 \begin{picture}(185,90)
  {
  \begin{picture}(180,45)(5,-45)
  {
 \begin{picture}(45,40)(5,0)    
 \put(5,5){\vector(1,0){45}}
 \put(5,5){\vector(0,1){35}}
 \put(20,35){\circle{2}}
 \put(28,35){\circle{2}}
 \put(45,35){\circle{2}}
 \put(20,35){\line(1,0){25}}
 \put(15,36){$X'$}
 \put(46,36){$X''$}
 \put(29,36){$Y$}
 \put(23.5,31){$\eps$}
  \end{picture}
 \begin{picture}(45,40)(-3,0)    
 \put(5,5){\vector(1,0){45}}
 \put(5,5){\vector(0,1){35}}
 \put(20,5){\circle{2}}
 \put(20,35){\circle{2}}
 \put(45,35){\circle{2}}
 \put(20,35){\line(1,0){25}}
 \put(20,5){\line(0,1){30}}
 \put(15,36){$X'$}
 \put(46,36){$X''$}
 \put(21,36){$Y''$}
 \put(22,6){$Y'$}
   \end{picture}
 \begin{picture}(25,40)(-11,0)    
 \put(5,5){\vector(1,0){25}}
 \put(5,5){\vector(0,1){35}}
 \put(5,6){\circle{2}}
 \put(20,6){\circle{2}}
 \put(20,35){\circle{2}}
 \put(20,6){\line(0,1){29}}
 \put(5,6){\line(1,0){15}}
 \put(-1,7){$X'$}
 \put(13,7){$X''$}
 \put(21,36){$Y''$}
 \put(22,6){$Y'$}
   \end{picture}
 \begin{picture}(25,40)(-19,0)     
 \put(5,5){\vector(1,0){25}}
 \put(5,5){\vector(0,1){35}}
 \put(20,6){\circle{2}}
 \put(5,6){\circle{2}}
 \put(5,6){\line(1,0){15}}
 \put(-1,7){$X'$}
 \put(15,7){$X''$}
 \put(6,7){$Y$}
   \end{picture}
 \begin{picture}(10,40)(-27,0)      
 \put(5,5){\vector(1,0){10}}
 \put(5,5){\vector(0,1){35}}
 \put(5,5){\circle{2}}
 \put(6,6){$Y$}
 \put(0,6){$X$}
   \end{picture}
   }
   \end{picture}
  \begin{picture}(180,45)(187,0)
  {
 \begin{picture}(45,40)(5,0)       
 \put(5,5){\vector(1,0){45}}
 \put(5,5){\vector(0,1){35}}
 \put(20,35){\circle{2}}
 \put(28,35){\circle{2}}
 \put(45,35){\circle{2}}
 \put(20,35){\line(1,0){25}}
 \put(15,36){$X'$}
 \put(46,36){$X''$}
 \put(29,36){$Y$}
 \put(23.5,31){$\eps$}
  \end{picture}
 \begin{picture}(45,40)(-3,0)        
 \put(5,5){\vector(1,0){45}}
 \put(5,5){\vector(0,1){35}}
 \put(5,35){\circle{2}}
 \put(28,35){\circle{2}}
 \put(36,35){\circle{2}}
 \put(5,35){\line(1,0){31}}
 \put(6.5,36){$X'$}
 \put(37.5,36){$X''$}
 \put(25,36){$Y$}
 \put(31,31){$\eps$}
   \end{picture}
 \begin{picture}(45,40)(-11,0)      
 \put(5,5){\vector(1,0){45}}
 \put(5,5){\vector(0,1){35}}
 \put(6,35){\circle{2}}
 \put(6,5){\circle{2}}
 \put(36,35){\circle{2}}
 \put(5,35){\line(1,0){31}}
 \put(6,5){\line(0,1){30}}
 \put(7.5,36){$X'$}
 \put(37.5,36){$X''$}
 \put(-1.5,34){$Y''$}
 \put(7,6){$Y'$}
   \end{picture}
 \begin{picture}(8,40)(-19,0)       
 \put(5,5){\vector(1,0){8}}
 \put(5,5){\vector(0,1){35}}
 \put(6,35){\circle{2}}
 \put(6,5){\circle{2}}
 \put(6,5){\line(0,1){30}}
 \put(0,6){$X$}
 \put(-1.5,34){$Y''$}
 \put(7,6){$Y'$}
   \end{picture}
 \begin{picture}(8,40)(-27,0)       
 \put(5,5){\vector(1,0){8}}
 \put(5,5){\vector(0,1){35}}
 \put(5,5){\circle{2}}
 \put(6,6){$Y$}
 \put(0,6){$X$}
   \end{picture}
   }
  \end{picture}
 }
   \end{picture}
  \end{center}
\caption{The transformations for an H-worm with $x''-x>x-x'$.
 In the upper line:  $v\to v_1\to v_2\to v_3\to v_4=O$.
 In the lower line: $v\to v'_1\to v'_2\to v'_3\to v'_4=O$.}
\label{fig:fig9}
  \end{figure}

 \medskip
{\bf Case 3}: $v$ is an H-worm with $x''-x>x-x'=:\eps$. This case is a bit
more complicated. The action of $\ttwo$ at $v$ moves $Y$ to $(x',y)$,
followed by moving $Y'$ to $(x',0)$. The action of $\tone$ at $v_1$ moves
$X''$ to $(x',y)$, followed by moving $X$ to $(x',0)$, followed by moving
$X'$ to $(0,0)$. The action of $\ttwo$ at $v_2$ moves $Y''$ to $(x',0)$,
followed by moving $Y$ to $(0,0)$. And the action of $\tone$ at $v_3$
moves $X''$ to $(0,0)$. (See the upper line of Fig.~\ref{fig:fig9}.)

In the second chain of transformations (see the lower line of
Fig.~\ref{fig:fig9}), the operator $\tone$ applied to $v$ moves $X''$ to
the point $(x+\eps,y)=(2x-x',y)$ (so that $Y$ becomes the mid-point of the
new interval $[X',X'']$). Then it moves $X'$ to $(0,y)$. The action of
$\ttwo$ at $v'_1$ moves $Y$ to $(0,y)$, followed by moving $Y'$ to
$(0,0)$. The action of $\tone$ at $v'_2$ moves $X''$ to $(0,y)$, followed
by moving $X$ to $(0,0)$. And the action of $\ttwo$ at $v'_3$ moves $Y''$
to $(0,0)$.

This gives $v_4=v'_4=O$ and:
   \begin{gather*}
  a_1=x-x'+y, \quad a_2=\frac{x''}{2}+y,\quad a_3=x'+y, \quad
  a_4=\frac{x'}{2}; \\
  b_1=\frac{x''-2\eps}{2}=\frac{x''}{2}+x'-x,\quad b_2=x+y,\quad
  b_3=\frac{x+\eps}{2}+y=x-\frac{x'}{2}+y, \quad b_4=y.
  \end {gather*}

Relation~\refeq{C1} turns into $x''+2y\ge x'+y\ge x'$, and~\refeq{C2}
turns into $x+y\ge x-\frac{x'}{2}+y\ge y$, which are valid.

Finally, we have $a_2\ge b_3$,
$p=(x''+2y)-(x'+y)+\frac{x'}{2}=x''-\frac{x'}{2}+y\ge b_3$, and
$q=(x-x'+y)+\frac{x'}{2}=x-\frac{x'}{2}+y=b_3$, yielding the first
equality in~\refeq{rel}. Also $a_1\ge b_4$, $r=(x''+2y)-(x'+y)=x''-x'+y\ge
b_4$, and $s=(x'+y)-x'=y=b_4$, yielding the second equality
in~\refeq{rel}.

 \medskip
{\bf Case 4}: $v$ is an H-worm with $2x\ge x'+x''$. The first chain of
transformations is similar to that in Case 3, giving similar
expressions for the numbers $a_i$. Compute the numbers $b_i$. The
action of $\tone$ at $v$ moves $X'$ to $(0,y)$. The action of $\ttwo$
at $v'_1$ moves $Y$ to $(0,y)$ and then moves $Y'$ to $(0,0)$. The
action of $\tone$ at $v'_2$ moves $X''$ to $(0,y)$ and then moves $X$
to $(0,0)$. And the action of $\ttwo$ at $v'_3$ moves $Y''$ to $(0,0)$.
This gives $v'_4=O$ and:
   $$
  b_1=\frac{x'}{2},\quad b_2=x+y,\quad
  b_3=\frac{x''}{2}+y, \quad b_4=y.
  $$

Relation~\refeq{C2} turns into $x+y\ge \frac{x''}{2}+y\ge y$, which is
valid since the condition $2x\ge x'+x''$ implies $x\ge \frac{x''}{2}$.

Finally, we have $a_2=b_3$, $p=x''-\frac{x'}{2}+y\ge b_3$, and
$q=x-\frac{x'}{2}+y\ge b_3$ (taking into account that $2x-x'\ge x''$),
yielding the first equality in~\refeq{rel}. Also $a_1\ge b_4$,
$r=x''-x'+y\ge b_4$, and $s=y=b_4$, yielding the second equality
in~\refeq{rel}.

 \medskip
{\bf Case 5}: $v$ is a V-worm with $2y\ge y'+y''$. By this condition,
the operator $\ttwo$ applied to $v$ moves $Y'$ to $(x,0)$. The action
of $\tone$ at $v_1$ moves $X$ to $(x,0)$ and then moves $X'$ to
$(0,0)$. The action of $\ttwo$ at $v_2$ moves $Y''$ to $(x,0)$ and then
moves $Y$ to $(0,0)$. And the action of $\tone$ at $v_3$ moves $X''$ to
$(0,0)$.

In the second chain, the action of $\tone$ at $v$ moves $X$ to $(x,y')$
($=Y'$) and then moves $X'$ to $(0,y')$. The action of $\ttwo$ at
$v'_1$ moves $Y''$ to $(x,y')$, then moves $Y$ to $(0,y')$, and
eventually moves $Y'$ to $(0,0)$. The action of $\tone$ at $v'_2$ moves
$X''$ to $(0,y')$ and then moves $X$ to $(0,0)$. And the action of
$\ttwo$ at $v'_3$ moves $Y''$ to $(0,0)$.

This gives $v_4=v'_4=O$ and:
   \begin{gather*}
  a_1=y', \quad a_2=\frac{x}{2}+y,\quad a_3=x+y'', \quad
  a_4=\frac{x}{2}; \\
  b_1=\frac{x}{2}+y-y',\quad b_2=x+y'',\quad b_3=\frac{x}{2}+y',
  \quad b_4=y'.
  \end {gather*}

Relation~\refeq{C1} turns into $x+2y\ge x+y''\ge x$, and~\refeq{C2}
turns into $x+y''\ge \frac{x}{2}+y'\ge y'$, which are valid (taking
into account that $2y\ge y''$).

Finally, we have $a_2\ge b_3$,
$p=(x+2y)-(x+y'')+\frac{x}{2}=\frac{x}{2}+2y-y''\ge b_3$ (since $2y-y''\ge
y'$), and $q=y'+\frac{x}{2}=b_3$, yielding the first equality
in~\refeq{rel}. Also $a_1=y'=b_4$, $r=(x+2y)-(x+y'')\ge 2y-y''\ge b_4$,
and $s=(x+y'')-x=y''\ge b_4$, yielding the second equality in~\refeq{rel}.

 \medskip
{\bf Case 6}: $v$ is a V-worm with $y''-y >y-y'=:\eps$. The action of
$\ttwo$ at $v$ moves $Y''$ to $(x,y+\eps)=(x,2y-y')$ (so that $X$
becomes the mid-point of the new interval $[Y',Y'']$). Then it moves
$Y'$ to $(x,0)$. The action of $\tone$ at $v_1$ moves $X$ to $(x,0)$
and then moves $X'$ to $(0,0)$. The action of $\ttwo$ at $v_2$ moves
$Y''$ to $(x,0)$ and then moves $Y$ to $(0,0)$. And the action of
$\tone$ at $v_3$ moves $X''$ to $(0,0)$. This gives $v_4=O$ and:
   $$
 a_1=y''-2\eps=y''-2y+2y', \quad a_2=\frac{x}{2}+y,\quad a_3=y+\eps+x=x+2y-y',
 \quad a_4=\frac{x}{2}.
   %
  $$

Relation~\refeq{C1} turns into $x+2y\ge x+2y-y'\ge x$, which is valid.

The second chain of transformations is similar to that in Case 5,
giving similar expressions for the numbers $b_i$ and validity
of~\refeq{C2}.

Finally, we have $a_2\ge \frac{x}{2}+y'=b_3$,
$p=(x+2y)-(x+2y-y')+\frac{x}{2}=y'+\frac{x}{2}=b_3$, and
$q=(y''-2y+2y')+\frac{x}{2}\ge y'+\frac{x}{2}=b_3$, yielding the first
equality in~\refeq{rel}. Also $a_1=y''-2y+2y'\ge y'=b_4$,
$r=(x+2y)-(x+2y-y')=y'=b_4$, and $s=(x+2y-y')-x\ge 2y-y'\ge y'=b_4$,
yielding the second equality in~\refeq{rel}.

 \medskip
 This completes the proof of Theorem~\ref{tm:W-B2}.

\section*{Appendix 2: Refining the ``local'' axioms }

In this section we translate the axioms for decorated $A_2$-crystals
given in Section~\ref{sec:axiom} so as to obtain implications in terms
of the original graph $G$. The resulting axioms (B5)--(B13) together
with their dual ones and Axioms (B0)--(B4),(B3$'$),(B4$'$), where (B1)--(B2) are replaced by their implementations (B1(i),(ii))--(B2(i),(ii)) will give a list of ``local'' axioms defining the regular (finite) $B_2$-crystals (which are applicable to their natural
infinite analogs though become not completely local). To make our description shorter, we do not give verbal formulations of the axioms, confining ourselves by merely
illustrating them in pictures, which to our belief is sufficient for
the reader to adequately restore the formal statements. To help, in each
illustration below we keep the same vertex notation as in the
illustration of the corresponding axiom in Section~\ref{sec:axiom},
with the only difference that for an $\otimes$-vertex $v$, the
corresponding central edge is denoted as $e_v$. Also in the left hand
sides of some pictures we do not indicate those edges whose existence
follow directly from Axioms (B3),(B3$'$),(B4),(B4$'$) and corollaries in
Section~\ref{sec:axiom}.

As before, the central vertices are indicated by thick dotes, and the
central edges by black rhombi.

The corresponding translations of parts (i),(ii) of (A0) have been
explained in Section~\ref{sec:axiom}. Also Axiom (A1) turns into
Corollary~\ref{cor:centr-edge}.

Axioms (B5)--(B13) below are derived from the {\em first} parts of
Axioms (A2)--(A8), whereas the second (symmetrical) parts of the latter
will derive the corresponding dual axioms (B5$'$)--(B13$'$).

The first part of (A2) turns into the following:

\unitlength=.8mm \special{em:linewidth 0.4pt} \linethickness{0.4pt}
\begin{picture}(140,35)(-20,0)
  \put(-20,20){{\bf (B5)}}
 \put(25,28){\circle{2}}
 \put(40,28){\circle{2}}
 \put(26,28){\vector(1,0){13}}
 \put(40,17){\vector(0,1){10}}
 \put(32,28){\makebox(0,0)[cc]{\scriptsize $\blacklozenge$}}
 \put(30,23){$e_v$}
 \put(63.00,18.00){\makebox(0,0)[cc]{$\Longrightarrow$}}
 \put(80,4){\circle{2}}
 \put(80,16){\circle{2}}
 \put(95,4){\circle{2}}
 \put(95,28){\circle{2}}
 \put(110,16){\circle{2}}
 \put(110,28){\circle{2}}
 \put(91,16){\circle*{2}}
 \put(95,20){\circle*{2}}
\put(81,4){\vector(1,0){13}}
 \put(81,16){\vector(1,0){9}}
 \put(92,16){\vector(1,0){17}}
 \put(96,28){\vector(1,0){13}}
 \put(80,5){\vector(0,1){10}}
 \put(95,5){\vector(0,1){14}}
 \put(95,21){\vector(0,1){6}}
 \put(110,17){\vector(0,1){10}}
  \put(87,4){\makebox(0,0)[cc]{\scriptsize $\blacklozenge$}}
  \put(102,28){\makebox(0,0)[cc]{\scriptsize $\blacklozenge$}}
 \put(100,23){$e_v$}
 \end{picture}

\smallskip
The first part of (A3), in the case when $e'$ leaves $v$, turns into
the following:

\unitlength=.8mm \special{em:linewidth 0.4pt} \linethickness{0.4pt}
\begin{picture}(140,32)(-20,0)
  \put(-20,20){{\bf (B6)}}
 \put(25,18){\circle{2}}
 \put(40,18){\circle{2}}
 \put(26,18){\vector(1,0){13}}
 \put(25,19){\vector(0,1){9}}
 \put(40,8){\vector(0,1){9}}
 \put(32,18){\makebox(0,0)[cc]{\scriptsize $\blacklozenge$}}
 \put(32,20){$e_v$}
 \put(63.00,18.00){\makebox(0,0)[cc]{$\Longrightarrow$}}
 \put(80,4){\circle*{2}}
 \put(80,14){\circle{2}}
 \put(80,18){\circle{2}}
 \put(80,28){\circle{2}}
 \put(95,4){\circle{2}}
 \put(95,18){\circle{2}}
 \put(100,14){\circle{2}}
 \put(100,28){\circle*{2}}
\put(81,4){\vector(1,0){13}}
 \put(81,14){\vector(1,0){18}}
 \put(81,18){\vector(1,0){13}}
 \put(81,28){\vector(1,0){18}}
 \put(80,5){\vector(0,1){8}}
 \put(80,19){\vector(0,1){8}}
 \put(95,5){\vector(0,1){12}}
 \put(100,15){\vector(0,1){12}}
  \put(90,14){\makebox(0,0)[cc]{\scriptsize $\blacklozenge$}}
  \put(87,18){\makebox(0,0)[cc]{\scriptsize $\blacklozenge$}}
 \put(87,20){$e_v$}
 \put(76,2){$u$}
 \end{picture}

\smallskip

The first part of (A3), in the case when $e'$ leaves $u$, turns into
the following:

\unitlength=.8mm \special{em:linewidth 0.4pt} \linethickness{0.4pt}
\begin{picture}(140,32)(-20,0)
  \put(-20,20){{\bf (B7)}}
 \put(25,18){\circle{2}}
 \put(40,18){\circle{2}}
 \put(25,4){\circle*{2}}
 \put(40,4){\circle{2}}
 \put(26,18){\vector(1,0){13}}
 \put(40,5){\vector(0,1){12}}
 \put(26,4){\vector(1,0){13}}
 \put(25,5){\vector(0,1){8}}
 \put(32,18){\makebox(0,0)[cc]{\scriptsize $\blacklozenge$}}
 \put(32,20){$e_v$}
 \put(20,2){$u$}
 \put(63.00,18.00){\makebox(0,0)[cc]{$\Longrightarrow$}}
 \put(80,4){\circle*{2}}
 \put(80,14){\circle{2}}
 \put(80,18){\circle{2}}
 \put(80,28){\circle{2}}
 \put(95,4){\circle{2}}
 \put(95,18){\circle{2}}
 \put(100,14){\circle{2}}
 \put(100,28){\circle*{2}}
\put(81,4){\vector(1,0){13}}
 \put(81,14){\vector(1,0){18}}
 \put(81,18){\vector(1,0){13}}
 \put(81,28){\vector(1,0){18}}
 \put(80,5){\vector(0,1){8}}
 \put(80,19){\vector(0,1){8}}
 \put(95,5){\vector(0,1){12}}
 \put(100,15){\vector(0,1){12}}
  \put(90,14){\makebox(0,0)[cc]{\scriptsize $\blacklozenge$}}
  \put(87,18){\makebox(0,0)[cc]{\scriptsize $\blacklozenge$}}
 \put(87,20){$e_v$}
 \put(75,2){$u$}
 \end{picture}

\smallskip
The first part of (A4) turns into the following:

\unitlength=.8mm \special{em:linewidth 0.4pt} \linethickness{0.4pt}
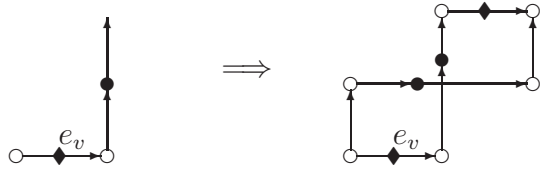
\begin{picture}(140,32)(-20,0)
  \put(-20,20){{\bf (B8)}}
 \put(25,4){\circle{2}}
 \put(40,4){\circle{2}}
 \put(40,16){\circle*{2}}
 \put(26,4){\vector(1,0){13}}
 \put(40,5){\vector(0,1){10}}
 \put(40,17){\vector(0,1){10}}
 \put(32,4){\makebox(0,0)[cc]{\scriptsize $\blacklozenge$}}
 \put(32,6){$e_v$}
 \put(42.5,15){$w$}
 \put(63.00,18.00){\makebox(0,0)[cc]{$\Longrightarrow$}}
 \put(80,4){\circle{2}}
 \put(80,16){\circle{2}}
 \put(95,4){\circle{2}}
 \put(95,28){\circle{2}}
 \put(110,16){\circle{2}}
 \put(110,28){\circle{2}}
 \put(91,16){\circle*{2}}
 \put(95,20){\circle*{2}}
\put(81,4){\vector(1,0){13}}
 \put(81,16){\vector(1,0){9}}
 \put(92,16){\vector(1,0){17}}
 \put(96,28){\vector(1,0){13}}
 \put(80,5){\vector(0,1){10}}
 \put(95,5){\vector(0,1){14}}
 \put(95,21){\vector(0,1){6}}
 \put(110,17){\vector(0,1){10}}
  \put(87,4){\makebox(0,0)[cc]{\scriptsize $\blacklozenge$}}
  \put(102,28){\makebox(0,0)[cc]{\scriptsize $\blacklozenge$}}
 \put(87,6){$e_v$}
 \end{picture}

\smallskip

The first part of (A5) turns into the following:

\unitlength=.8mm \special{em:linewidth 0.4pt} \linethickness{0.4pt}
\begin{picture}(140,32)(-20,0)
  \put(-20,20){{\bf (B9)}}
 \put(25,4){\circle*{2}}
 \put(40,4){\circle{2}}
 \put(40,16){\circle*{2}}
 \put(26,4){\vector(1,0){13}}
 \put(40,5){\vector(0,1){10}}
 \put(20,2){$u$}
 \put(42,17){$v$}
 \put(63.00,18.00){\makebox(0,0)[cc]{$\Longrightarrow$}}
 \put(80,4){\circle*{2}}
 \put(95,4){\circle{2}}
 \put(95,16){\circle*{2}}
 \put(95,28){\circle{2}}
 \put(110,28){\circle*{2}}
 \put(81,4){\vector(1,0){13}}
 \put(96,28){\vector(1,0){13}}
 \put(95,5){\vector(0,1){10}}
 \put(95,17){\vector(0,1){10}}
 \put(75,2){$u$}
 \put(97,17){$v$}
 \end{picture}

\smallskip
The first part of (A6), in the case when $e'$ enters $u$, turns into the following:

\unitlength=.8mm \special{em:linewidth 0.4pt} \linethickness{0.4pt}
\begin{picture}(140,35)(0,0)
  \put(0,20){{\bf (B10)}}
 \put(45,16){\circle*{2}}
 \put(60,16){\circle{2}}
 \put(60,28){\circle*{2}}
 \put(31,16){\vector(1,0){13}}
 \put(46,16){\vector(1,0){13}}
 \put(60,17){\vector(0,1){10}}
 \put(42,12){$u$}
 \put(62,29){$v$}
 \put(83.00,18.00){\makebox(0,0)[cc]{$\Longrightarrow$}}
 \put(100,4){\circle*{2}}
 \put(100,16){\circle{2}}
 \put(115,4){\circle{2}}
 \put(115,28){\circle{2}}
 \put(130,16){\circle{2}}
 \put(130,28){\circle*{2}}
 \put(111,16){\circle*{2}}
 \put(115,20){\circle*{2}}
\put(101,4){\vector(1,0){13}}
 \put(101,16){\vector(1,0){9}}
 \put(112,16){\vector(1,0){17}}
 \put(116,28){\vector(1,0){13}}
 \put(100,5){\vector(0,1){10}}
 \put(115,5){\vector(0,1){14}}
 \put(115,21){\vector(0,1){6}}
 \put(130,17){\vector(0,1){10}}
 \put(108,12){$u$}
 \put(132,29){$v$}
 \end{picture}

\smallskip
The first part of (A6), in the case when $e'$ enters $v$, turns into the following:

\unitlength=.8mm \special{em:linewidth 0.4pt} \linethickness{0.4pt}
\begin{picture}(140,35)(0,0)
  \put(0,20){{\bf (B11)}}
 \put(45,16){\circle*{2}}
 \put(60,16){\circle{2}}
 \put(60,28){\circle*{2}}
 \put(46,16){\vector(1,0){13}}
 \put(46,28){\vector(1,0){13}}
 \put(60,17){\vector(0,1){10}}
 \put(42,12){$u$}
 \put(62,29){$v$}
 \put(83.00,18.00){\makebox(0,0)[cc]{$\Longrightarrow$}}
 \put(100,4){\circle*{2}}
 \put(100,16){\circle{2}}
 \put(115,4){\circle{2}}
 \put(115,28){\circle{2}}
 \put(130,16){\circle{2}}
 \put(130,28){\circle*{2}}
 \put(111,16){\circle*{2}}
 \put(115,20){\circle*{2}}
\put(101,4){\vector(1,0){13}}
 \put(101,16){\vector(1,0){9}}
 \put(112,16){\vector(1,0){17}}
 \put(116,28){\vector(1,0){13}}
 \put(100,5){\vector(0,1){10}}
 \put(115,5){\vector(0,1){14}}
 \put(115,21){\vector(0,1){6}}
 \put(130,17){\vector(0,1){10}}
 \put(108,12){$u$}
 \put(132,29){$v$}
 \end{picture}

\smallskip
The first part of (A7) turns into the following:

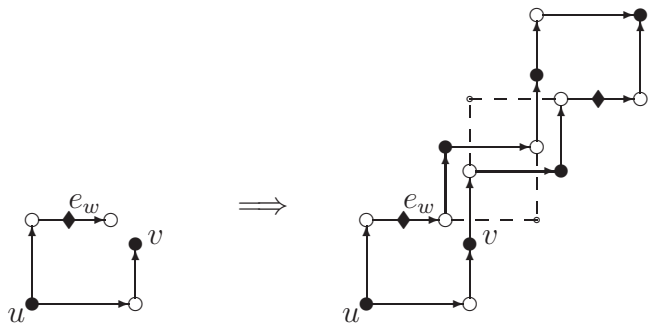
\unitlength=.8mm \special{em:linewidth 0.4pt} \linethickness{0.4pt}
\begin{picture}(150,57)(0,0)
  \put(0,40){{\bf (B12)}}
 \put(45,4){\circle*{2}}
 \put(45,18){\circle{2}}
 \put(58,18){\circle{2}}
 \put(62,4){\circle{2}}
 \put(62,14){\circle*{2}}
 \put(46,4){\vector(1,0){15}}
 \put(46,18){\vector(1,0){11}}
 \put(45,5){\vector(0,1){12}}
 \put(62,5){\vector(0,1){8}}
  \put(51,18){\makebox(0,0)[cc]{\scriptsize $\blacklozenge$}}
 \put(41,1){$u$}
 \put(64,14){$v$}
 \put(51,20){$e_w$}
 \put(83,20){\makebox(0,0)[cc]{$\Longrightarrow$}}
 \put(100,4){\circle*{2}}
 \put(100,18){\circle{2}}
 \put(113,18){\circle{2}}
 \put(113,30){\circle*{2}}
 \put(117,4){\circle{2}}
 \put(117,14){\circle*{2}}
 \put(117,26){\circle{2}}
 \put(128,30){\circle{2}}
 \put(128,42){\circle*{2}}
 \put(128,52){\circle{2}}
 \put(132,26){\circle*{2}}
 \put(132,38){\circle{2}}
 \put(145,38){\circle{2}}
 \put(145,52){\circle*{2}}
 \put(101,4){\vector(1,0){15}}
 \put(101,18){\vector(1,0){11}}
 \put(118,26){\vector(1,0){13}}
 \put(114,30){\vector(1,0){13}}
 \put(133,38){\vector(1,0){11}}
 \put(129,52){\vector(1,0){15}}
 \put(100,5){\vector(0,1){12}}
 \put(113,19){\vector(0,1){10}}
 \put(117,5){\vector(0,1){8}}
 \put(117,15){\vector(0,1){10}}
  \put(128,31){\vector(0,1){10}}
 \put(128,43){\vector(0,1){8}}
 \put(132,27){\vector(0,1){10}}
 \put(145,39){\vector(0,1){12}}
  \put(106,18){\makebox(0,0)[cc]{\scriptsize $\blacklozenge$}}
  \put(138,38){\makebox(0,0)[cc]{\scriptsize $\blacklozenge$}}
 \put(117,38){\circle{1}}
 \put(128,18){\circle{1}}
 \put(117,27){\line(0,1){2}}
 \put(117,31){\line(0,1){2}}
 \put(117,35){\line(0,1){2}}
 \put(128,19){\line(0,1){2}}
 \put(128,23){\line(0,1){2}}
 \put(128,27){\line(0,1){2}}
 \put(114,18){\line(1,0){2}}
 \put(118,18){\line(1,0){2}}
 \put(122,18){\line(1,0){2}}
 \put(126,18){\line(1,0){2}}
 \put(117,38){\line(1,0){2}}
 \put(121,38){\line(1,0){2}}
 \put(125,38){\line(1,0){2}}
 \put(129,38){\line(1,0){2}}
 \put(96,1){$u$}
 \put(119,14){$v$}
 \put(106,20){$e_w$}
   \end{picture}

\smallskip
\noindent (Here the required configuration in the right hand side is
drawn by solid lines, while the dashed lines indicate the
fragments that are automatically added by the commutativity axioms
(B3),(B3$'$). This gives the crystal $S(1,1)$.)

\smallskip
The first part of (A8) turns into the following:

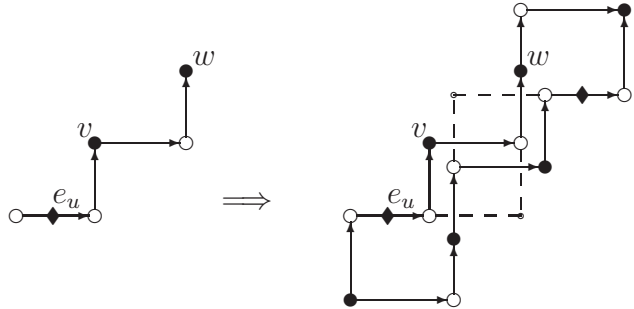
\unitlength=.8mm \special{em:linewidth 0.4pt} \linethickness{0.4pt}
\begin{picture}(150,57)(0,0)
  \put(0,40){{\bf (B13)}}
 \put(45,18){\circle{2}}
 \put(58,18){\circle{2}}
 \put(58,30){\circle*{2}}
 \put(73,30){\circle{2}}
 \put(73,42){\circle*{2}}
 \put(46,18){\vector(1,0){11}}
 \put(58,19){\vector(0,1){10}}
 \put(59,30){\vector(1,0){13}}
  \put(73,31){\vector(0,1){10}}
  \put(51,18){\makebox(0,0)[cc]{\scriptsize $\blacklozenge$}}
 \put(51,20){$e_u$}
 \put(55,31){$v$}
 \put(74,43){$w$}
 \put(83,20){\makebox(0,0)[cc]{$\Longrightarrow$}}
 \put(100,4){\circle*{2}}
 \put(100,18){\circle{2}}
 \put(113,18){\circle{2}}
 \put(113,30){\circle*{2}}
 \put(117,4){\circle{2}}
 \put(117,14){\circle*{2}}
 \put(117,26){\circle{2}}
 \put(128,30){\circle{2}}
 \put(128,42){\circle*{2}}
 \put(128,52){\circle{2}}
 \put(132,26){\circle*{2}}
 \put(132,38){\circle{2}}
 \put(145,38){\circle{2}}
 \put(145,52){\circle*{2}}
 \put(101,4){\vector(1,0){15}}
 \put(101,18){\vector(1,0){11}}
 \put(118,26){\vector(1,0){13}}
 \put(114,30){\vector(1,0){13}}
 \put(133,38){\vector(1,0){11}}
 \put(129,52){\vector(1,0){15}}
 \put(100,5){\vector(0,1){12}}
 \put(113,19){\vector(0,1){10}}
 \put(117,5){\vector(0,1){8}}
 \put(117,15){\vector(0,1){10}}
  \put(128,31){\vector(0,1){10}}
 \put(128,43){\vector(0,1){8}}
 \put(132,27){\vector(0,1){10}}
 \put(145,39){\vector(0,1){12}}
  \put(106,18){\makebox(0,0)[cc]{\scriptsize $\blacklozenge$}}
  \put(138,38){\makebox(0,0)[cc]{\scriptsize $\blacklozenge$}}
 \put(117,38){\circle{1}}
 \put(128,18){\circle{1}}
 \put(117,27){\line(0,1){2}}
 \put(117,31){\line(0,1){2}}
 \put(117,35){\line(0,1){2}}
 \put(128,19){\line(0,1){2}}
 \put(128,23){\line(0,1){2}}
 \put(128,27){\line(0,1){2}}
 \put(114,18){\line(1,0){2}}
 \put(118,18){\line(1,0){2}}
 \put(122,18){\line(1,0){2}}
 \put(126,18){\line(1,0){2}}
 \put(117,38){\line(1,0){2}}
 \put(121,38){\line(1,0){2}}
 \put(125,38){\line(1,0){2}}
 \put(129,38){\line(1,0){2}}
 \put(106,20){$e_u$}
 \put(110,31){$v$}
 \put(129,43){$w$}
   \end{picture}

\medskip
It should be noted that Axioms (B0)--(B13),(B3$'$)--(B13$'$) characterize a
slightly larger class than the class of (finite and infinite) S-graphs
because we did not translate parts (iii),(iv) of (A0) (their
implementations in local terms hardly exist). The arising extra graphs
are infinite and may be named as ``exotic'' infinite analogs of regular
$B_2$-crystals. Such graphs are obtained by the construction in
Section~\ref{sec:constr} when we take as $C$ one of the exotic infinite
analogs of $A_2$-crystals (having no principal points at all),
described in~\cite[Sect.~6]{A-2}, and replace the monochromatic strings in it by
``$B$-sails'' with constant slopes. An instance is the ``left
$B$-sail'' with slope 2 which is unbounded from below and bounded from
above.


\end{document}